\documentclass[11pt]{article}
\usepackage[margin=1in]{geometry}
\usepackage{amsmath,amsfonts,amssymb,amsthm,mathrsfs}
\usepackage{graphicx}
\usepackage{subcaption}
\usepackage{bbm}
\usepackage{booktabs}
\usepackage{authblk}
\usepackage{hyperref}

\newcommand{\norm}[1]{\left\lVert#1\right\rVert}

\theoremstyle{plain}
\newtheorem{theorem}{Theorem}[section]
\newtheorem{lemma}[theorem]{Lemma}
\newtheorem{proposition}[theorem]{Proposition}
\newtheorem{corollary}[theorem]{Corollary}

\theoremstyle{definition}

\newtheorem{assumption}{Assumption}
\newtheorem{remark}[theorem]{Remark}

\title{Metropolis Adjusted Implicit Midpoint Langevin Algorithm: A fully implicit ergodic MCMC method}
\author[1,2]{P. Dobson}
\author[2,6]{A. Harrison}
\author[3]{T. Klatzer}
\author[4]{G. Roberts}
\author[5]{J. M. Sanz-Serna}
\author[2,6]{K. C. Zygalakis}
\affil[1]{\small School of Mathematics and Computer Science, Heriot-Watt University, EH14 4AS, Edinburgh, Scotland, UK}
\affil[2]{\small Maxwell Institute for Mathematical Sciences, Bayes Centre, 47 Potterrow, Edinburgh, Scotland, UK}
\affil[3]{\small School of Mathematical Sciences, Lancaster University, Lancaster, UK}
\affil[4]{\small Department of Statistics, University of Warwick, Coventry, CV4 7AL}
\affil[5]{\small Departamento de Matem\'aticas, Universidad Carlos III de Madrid, Avenida Universidad 30, 28911 Legan\'es, Madrid}
\affil[6]{\small School of Mathematics,  University of Edinburgh, Peter Guthrie Tait Rd, EH9 3FD, Edinburgh, Scotland, UK}
\date{}

\begin{document}
\maketitle

\begin{abstract}
We study the implicit midpoint Langevin algorithm (IMLA) and its Metropolis-adjusted variant (MAIMLA) for sampling probability distributions, possibly with non-Gaussian tails. The Metropolis--Hastings acceptance rate formula for the implicit midpoint proposal does not require the computation of determinants. We establish geometric ergodicity of IMLA under dissipativity and polynomial growth conditions on the gradient; the hypotheses accommodate potentials with superlinearly growing gradients. For MAIMLA, we prove geometric ergodicity for convex potentials in two settings. The first corresponds to superlinear gradients with radially symmetric leading-order behaviour; this includes targets for which MALA is not geometrically ergodic. The second deals with globally Lipschitz gradients with a suitable decay condition on the third derivatives of the potential.
 We also derive high-dimensional optimal-scaling limits for MAIMLA on targets that are products of independent copies, both for identically distributed and scaled factors.  A number of numerical experiments illustrate the  improved sampling efficiency of MAIMLA over MALA in some settings.

\medskip
\noindent\textbf{MSC2020 subject classifications:} Primary 65C05, 60J22; secondary 65C40, 60H35, 65L04.

\end{abstract}

\section{Introduction}
Sampling from a target distribution of the form
\begin{equation*}
\pi(x)\propto \exp(-U(x)), \qquad x\in\mathbb{R}^{d},
\end{equation*}
is a fundamental problem in computational statistics \cite{BPRSSS13,RR98,RT96}, Bayesian inverse problems \cite{S10}, molecular dynamics \cite{LM15}, and uncertainty quantification \cite{BLSZ18,Smith13}. A common approach to deal with the sampling problem is to construct a Markov process on an extended state space $\mathbb{R}^{m}$, with $m\geq d$, whose invariant distribution, or an appropriate marginal thereof, coincides with the target distribution $\pi$.

A starting point for the derivation of many MCMC algorithms is the overdamped Langevin equation
\begin{equation} \label{eq:Langevin}
dX_{t} = -\nabla U(X_{t})\,dt + \sqrt{2}\,dW_{t},
\end{equation}
where $X_{t} \in \mathbb{R}^{d}$ and $W_{t}$ is the standard $d$-dimensional Brownian motion. This equation has the property \cite{RT96}
that, under mild conditions on the potential $U$, the law of $X_t$ converges to $\pi$
as $t\to\infty$, and does so exponentially fast (see Theorem~\ref{thm:langevin_geometric_ergodicity} below for a precise statement). With the
exception of a few special cases, however, \eqref{eq:Langevin} cannot be
solved exactly, and in practice one resorts to numerical simulation. This
introduces two separate difficulties. First, the discretisation must
inherit the stability of the continuous dynamics; second, the discrete
chain has an invariant measure that differs from $\pi$, and this bias
must either be controlled or removed. The standard way to remove it is a
Metropolis--Hastings correction, which turns the discretisation into a
proposal for a Markov chain that is exactly $\pi$-invariant; the
Euler--Maruyama discretisation corrected in this way gives the Metropolis
Adjusted Langevin Algorithm (MALA) \cite{RT96}.

A typical assumption when establishing geometric ergodicity both for unadjusted schemes \cite{MSH02, DM17}
and for MALA \cite{RT96, DM22},
is that $\nabla U$ is globally Lipschitz, so that $\pi$ has tails that
are at most Gaussian. Many targets of interest violate this requirement and unfortunately explicit integrators run into difficulties whenever $\pi$ is \emph{light tailed} and
$\nabla U$ grows superlinearly.  On a bounded
time interval, the Euler--Maruyama scheme may diverge in both the strong
and the weak sense \cite{HJK11}; the resulting Markov chain may be transient
and its Metropolised version MALA
may fail to be geometrically ergodic for every  fixed step size
\cite{RT96}.
The one-dimensional family $\pi(x)\propto\exp(-c|x|^{\beta})$, $\beta>2$, $c>0$,
already exhibits all of these pathologies.

Many different remedies have been proposed to fix the issue of geometric ergodicity in the superlinear growth case for explicit integrators. \emph{Taming} modifies the drift so
that it remains bounded, or grows at most linearly, which restores
convergence of the unadjusted scheme \cite{HJK12, S13, BDMS19}
and geometric ergodicity of the adjusted one \cite{Atchade06};  more recently, the Barker proposal
of \cite{LZ22} achieves robustness to tail behaviour by
replacing the Gaussian proposal by a skew-symmetric one. A different
route is to make the integrator \emph{implicit}
\cite{KP92, H00}. Implicit schemes are indeed stable for light-tailed
targets \cite{CRS11},
but their Metropolisation may be problematic: the probability density of an implicit
proposal contains a Jacobian determinant, so that the
acceptance ratio in principle involves a ratio of such determinants
\cite{HSR21},
which  is likely to be prohibitively expensive to compute in high dimension. The partially
 linearly solvable proposals of
\cite{CRS11}
avoid this and are geometrically ergodic for the one-dimensional
super-Gaussian family, but this is achieved by treating implicitly only part of the drift. Related ideas appear in proximal MCMC
\cite{P16, KDAPSZ24},
where an implicit (proximal) step is applied to a possibly nonsmooth
part of $U$, and also in the Crank--Nicolson proposals for targets that are
changes of measure from a Gaussian reference distribution; these proposals are exact (and therefore would always be accepted) in the case where the target reduces to   the Gaussian reference
\cite{BRS08, CRSW13}.

In this paper, we study the geometric ergodicity of the implicit midpoint Langevin algorithm (IMLA) \cite{KDAPSZ24} and its Metropolis-adjusted variant (MAIMLA). Two structural features of
the midpoint rule make it attractive as a proposal. The first one is that, similarly to the situation for the unadjusted Langevin algorithm (ULA), IMLA can be interpreted as a one step discretisation of an underlying Hamiltonian system. When so  interpreted, the midpoint integrator is reversible and symplectic, so that the
standard Hamiltonian Monte Carlo formula for the acceptance rate \cite{BRSS18}
applies, and this entails that the computation of the acceptance rate does not involve any determinants.
The reason behind this  is  the volume preservation of the
symplectic map, which removes the Jacobian that typically burdens  implicit schemes. A second attractive feature is that  it is well known that, when $U$ is quadratic, the implicit midpoint exactly preserves the value of the Hamiltonian, and this implies that all proposals are accepted,
in the same way as
Crank--Nicolson proposals are. As a result, high acceptance rates may be expected in parts of the state space where the target is approximately Gaussian.

Our main results are the following.
\begin{enumerate}
\item \emph{Geometric ergodicity of IMLA} (Theorem~\ref{thm:geometric_ergodicity_imla}) for every step size $\delta>0$ for
convex $U$ satisfying a dissipativity condition and a polynomial growth
bound on $\nabla U$. The growth exponent
may be arbitrary, so superlinear gradients are covered.
\item \emph{Geometric ergodicity of MAIMLA in the superlinear regime}
(Theorem~\ref{thm:geometric_ergodicity_MAIMLA}) for every step size $\delta>0$ for convex potentials whose leading-order
term is radially symmetric: $U(x)=C\norm{x}^{1+p}+U_0(x)$ ($p>1$ and
$\nabla U_0$ of lower order). In particular,
for $\pi(x)\propto\exp(-c|x|^{\beta})$ with $\beta>2$, $c>0$, MAIMLA is
geometrically ergodic for all step sizes, whereas MALA is geometrically
ergodic for none \cite{RT96}.
We show by means of an example that the symmetry assumption cannot simply be dropped
(Remark~\ref{rem:asym}).
\item \emph{Geometric ergodicity of MAIMLA in the Lipschitz regime}
(Theorem~\ref{thm:maimla-lip}), for convex $U$ with $L$-Lipschitz
gradient, provided that $\nabla^3U$ decays at least cubically at infinity and
$\delta<2/L$. Since the acceptance probability depends on $U$ only
through its third derivatives, some condition of this type is
unavoidable; we exhibit a strongly log-concave target with Lipschitz
gradient, outside this class, for which the acceptance probability of
MAIMLA decays to zero in the tails.
\item \emph{Optimal scaling}. Using the
Hamiltonian interpretation and the framework of \cite{OptimalScalingSymmetry26},
we derive the high-dimensional limits of the acceptance probability and
of the expected squared jumping distance for MAIMLA on identically distributed and
scaled independent product targets.  In the i.i.d.\ case the step size scales as $d^{-1/3}$ as $d\uparrow\infty$ and the optimal acceptance probability is
$\approx 0.574$, just as in MALA.   However, the key integrator-dependent constant in the analysis is always smaller for MAIMLA than for MALA   and this implies that MAIMLA requires smaller step sizes than MALA to achieve a target acceptance probability and that, for a given acceptance probability, MAIMLA achieves larger expected jumping distances.
\end{enumerate}

The rest of the paper is organised as follows. Section \ref{sec:Langevin} reviews the overdamped Langevin equation,
its explicit and $\theta$-implicit discretisations, and  known
geometric ergodicity results for MALA. Section \ref{sec:maimla} introduces IMLA and MAIMLA, establishes the Hamiltonian interpretation
and well-posedness of the implicit step, and derives the optimal-scaling
limits. Section~\ref{sec:geo} presents geometric ergodicity results first for IMLA and then for MAIMLA in the superlinear and
Lipschitz regimes. Section~\ref{sec:num} is devoted to  numerical
experiments and Section~\ref{sec:con} concludes.  Bounds on the
solution of the midpoint equation, used throughout
Section~\ref{sec:geo}, are collected in Appendix~\ref{sec:midpoint}.

\section{Overdamped Langevin equation and its numerical discretisations} \label{sec:Langevin}

The different geometric ergodicity results in the paper require
 different assumptions on $U(x)$ and  $\nabla U(x)$ and for convenience we group these assumptions here.
\begin{assumption} \label{ass:ass1}
The function $U \in \mathcal{C}^{2}$ satisfies
\begin{equation} \label{eq:strong_ball}
\left \langle \nabla U(x),x \right \rangle   \geq a \norm{x}^{2}+b, \quad \forall x \in \mathbb{R}^{d}
\end{equation}
for some $a>0, b\in \mathbb{R}$.
\end{assumption}

\begin{assumption} \label{ass:ass2}
The function $U \in \mathcal{C}^{2}$ satisfies
 \begin{equation}\label{eq:V_growth_cond}
    \norm{\nabla U(x)} \leq c_p \norm{x}^p+C_p, \quad \forall x \in \mathbb{R}^{d}
\end{equation}
for some $p\geq1, c_p>0, C_p\in\mathbb{R}$.
\end{assumption}

\begin{assumption} \label{ass:ass3}
$U \in \mathcal{C}^2$ can be written as
    \begin{equation}\label{eq:radial_condition}
        U(x) = C \norm{x}^{1+p}+U_0(x)
    \end{equation}
     with $C>0$, $p>1$ and $U_0$ satisfying $\norm{\nabla U_0(x)}\leq C_0\norm{x}^{p-1}+c_0$, $C_0,c_0 \geq 0$.
     Note that $\nabla U(x)$ grows superlinearly with $x$.
\end{assumption}

\begin{assumption} \label{ass:ass4}
The function $U\in\mathcal{C}^2(\mathbb{R}^d)$ satisfies the following conditions.
There exist constants $L>0$, $m>0$, and $R>0$ such that
\begin{equation}
\norm{\nabla U(x)-\nabla U(y)} \leq L \norm{x-y},
\qquad x,y\in\mathbb{R}^d,
\end{equation}
and $U$ is strongly convex outside the ball of radius $R$, in the sense that
\begin{equation}
v^{\top}\nabla^2 U(x)v \geq m\norm{v}^2,
\qquad \norm{x}\geq R,\quad v\in\mathbb{R}^d,
\end{equation}
and $\nabla U(0)=0$.
\end{assumption}

\begin{assumption}\label{ass:third_derivative_decay}
The function $U\in \mathcal{C}^3(\mathbb{R}^d)$  has an $L$-Lipschitz gradient, and satisfies
\begin{equation}\label{eq:third_decay}
    \sup_{x\in\mathbb{R}^d}\norm{\nabla^3 U(x)} (1+\norm{x})^3 \leq C.
\end{equation}
\end{assumption}

We now review some basic facts about the Langevin equation \eqref{eq:Langevin}
and its numerical discretisations.

For the continuous-time Langevin dynamics we have the following result \cite[Theorem~2.2]{RT96}. Here and later, for a measurable $V:\mathbb{R}^d\to[1,\infty)$ and a signed measure $\mu$ we set
$\norm{\mu}_V=\sup_{|f|\le V}|\mu(f)|$.

\begin{theorem}[Geometric ergodicity of overdamped Langevin dynamics]
\label{thm:langevin_geometric_ergodicity}
Suppose that $U \in \mathcal{C}^{2}(\mathbb{R}^{d})$ satisfies
Assumption~\ref{ass:ass1}. Then,
 there exist constants $K<\infty$ and
$\rho>0$ such that
\[
    \norm{P_t(x,\cdot)-\pi}_V
    \leq K V(x)\exp(-\rho t),
    \qquad x\in\mathbb{R}^d,\quad t\geq0,
\]
where  $P_t(x,\cdot)$ is the transition kernel of \eqref{eq:Langevin}, $\pi$ is the target, and $V(x)=1+\norm{x}^{2}.$
\end{theorem}
Consequently, the distribution of the continuous dynamics approaches the target exponentially fast. Since exact simulation is generally unavailable, we  consider numerical discretisations. These introduce two distinct difficulties: preserving the stability of the continuous dynamics and controlling the bias in the invariant distribution. The simplest numerical method for solving \eqref{eq:Langevin} is the Euler--Maruyama method
\begin{equation} \label{eq:euler_maruyama}
X_{n+1}=X_{n}-\delta \nabla U(X_{n})+\sqrt{2\delta}\xi_{n},
\end{equation}
where $(\xi_n)_{n\geq 0}$ is a sequence of i.i.d.\ $\mathcal{N}(0, I_d)$ random variables. A natural question is whether this discretisation preserves the invariant distribution $\pi$. In general, it does not \cite{KZ11}: even for a positive definite quadratic potential, the Euler--Maruyama chain converges, for $\delta$ sufficiently small, to an invariant density $\pi^{\delta}$ that differs from $\pi$.

One way to remove this invariant-measure bias is to add a Metropolis--Hastings correction, which preserves the target distribution $\pi$. In particular, given the current state $X_n=x$, propose
\begin{equation*}
    y=x-\delta\nabla U(x)+\sqrt{2\delta}\,\xi_n
\end{equation*}
and accept the proposal with probability
\begin{align*}
    \alpha_\delta(x,y)
    &=\min\left\{1,\frac{\pi(y)q_\delta(y,x)}
                            {\pi(x)q_\delta(x,y)}\right\},
\end{align*}
where $q_{\delta}(x,y)$ is the density at $y$ of $\mathcal{N}(x-\delta\nabla U(x),2\delta I_d)$. This algorithm is known as the Metropolis Adjusted Langevin Algorithm (MALA). The following geometric ergodicity result follows from \cite[Theorem~1]{DM22}.
\begin{theorem}[Geometric ergodicity of MALA]
\label{thm:mala_geometric_ergodicity}
Suppose that $U$ satisfies Assumption~\ref{ass:ass4}.
Then there exists $\delta_0>0$ such that, for every
$\delta\in(0,\delta_0]$, there exist constants
$K_\delta<\infty$ and $\rho_\delta\in(0,1)$ such that
\[
    \norm{(P^{\text{MALA}}_\delta)^n(x,\cdot)-\pi}_V
    \leq K_\delta V(x)\rho_\delta^n,
    \qquad x\in\mathbb{R}^d,\quad n\geq0,
\]
where $P^{\text{MALA}}_\delta$ is the transition kernel of MALA
with step size $\delta$, $\pi$ is the target distribution,
and
\[
    V(x)=\exp\left(\frac{m}{16}\norm{x}^2\right),
\]
with $m$ the strong convexity constant in
Assumption~\ref{ass:ass4}.
\end{theorem}

Theorem~\ref{thm:mala_geometric_ergodicity} requires considerably more demanding hypotheses than the corresponding result for the
continuous dynamics. In particular, its global Lipschitz condition forces $U$ to grow at most quadratically, so that $\pi$ has tails that
are at least as heavy as Gaussian; potentials whose gradient grows
superlinearly, and hence targets with lighter-than-Gaussian tails, are
excluded. This exclusion is not an artefact of the proof. Consider the
one-dimensional family
\begin{equation}\label{eq:supergaussian}
\pi(x)\propto\exp\big(-c\,|x|^{\beta}\big),\qquad c>0,\ \beta>2 ,
\end{equation}
for which $U'(x)=c\beta\,\mathrm{sign}(x)|x|^{\beta-1}$ grows
superlinearly.  For these targets, the Euler--Maruyama chain \eqref{eq:euler_maruyama}
is transient \cite[Theorem~3.2(b)]{RT96},
and MALA is not geometrically ergodic, regardless of how small the fixed
step size $\delta$ is chosen \cite[Theorem~4.2]{RT96}.  Furthermore, even on a bounded time interval,
Euler--Maruyama approximations can diverge in the strong $L^p$
and  weak senses as the step size tends to zero
\cite{HJK11}.

A natural approach to improving stability when $\nabla U$ grows superlinearly is to consider implicit numerical discretisations of \eqref{eq:Langevin}. One such family is the stochastic $\theta$-scheme \cite{H00,HSR21},
\begin{equation} \label{eq:standard_theta}
X_{n+1}=X_{n}-\delta (1-\theta)\nabla U(X_{n})-\delta \theta \nabla U(X_{n+1})+\sqrt{2\delta}\xi_{n},
\end{equation}
where $\theta \in [0,1]$. Implicitness can restore geometric
ergodicity in settings where the explicit Euler scheme fails.
For example, for targets of the form \eqref{eq:supergaussian}
the unadjusted scheme is geometrically ergodic for
$\theta>1/2$ and any $\delta>0$ \cite[Section~3.2]{CRS11}.
Its invariant distribution, however, generally differs
from the target $\pi$.

A Metropolis--Hastings correction can remove this bias,
but introduces an additional computational cost. Provided
the implicit proposal is globally well-defined through an
invertible map, its density contains the
Jacobian factor
$\left|\det(I_d+\delta\theta\nabla^2 U(y))\right|$
\cite{HSR21} and, as a consequence, the Metropolis--Hastings formula for the acceptance
ratio involves
\[
    \frac{\left|\det(I_d+\delta\theta\nabla^2 U(x))\right|}
         {\left|\det(I_d+\delta\theta\nabla^2 U(y))\right|}.
\]
Thus, even if geometric ergodicity of the adjusted scheme
were known, evaluating these Hessian determinants
could make the algorithm computationally unattractive in high dimensions. An alternative is to use partially implicit, linearly
solvable proposals. Such methods are developed in \cite{CRS11},
where geometric ergodicity is established for a Metropolis-adjusted
split-step variant targeting the one-dimensional super-Gaussian
densities above when $\theta>1/2$. The implicit midpoint rule studied next retains full implicitness in the drift while avoiding determinants entirely.

\section{Metropolis Adjusted Implicit Midpoint Langevin  Algorithm}
\label{sec:maimla}

In this section we introduce the implicit midpoint discretisation of
\eqref{eq:Langevin} and its Metropolis adjustment, and we establish the
structural properties that will be used throughout the paper.

\subsection{Definition and well-posedness}

The \emph{implicit midpoint Langevin algorithm} (IMLA)
\cite{KDAPSZ24}
is the discretisation
\begin{equation}\label{eq:imla}
X_{n+1}=X_n-\delta\,\nabla U\!\Big(\frac{X_{n+1}+X_n}{2}\Big)+\sqrt{2\delta}\,\xi_n ,
\qquad n\ge 0,
\end{equation}
of \eqref{eq:Langevin}. For quadratic potentials \eqref{eq:imla} coincides
with the stochastic $\theta$-scheme \eqref{eq:standard_theta} with $\theta=1/2$;
for general potentials the two schemes differ, and it is the midpoint
form of the drift that will allow a determinant-free Metropolis
correction below.

It is convenient, both for the implementation and the analysis, to write one step of \eqref{eq:imla} in terms of the
midpoint value. Given the current state $x$ and noise $\xi$, let
$x_{1/2}$ denote the solution of the \emph{midpoint equation}
\begin{equation}\label{eq:midpoint}
x_{1/2}+\frac{\delta}{2}\,\nabla U(x_{1/2})=x+\sqrt{\frac{\delta}{2}}\;\xi ,
\end{equation}
and set $x^\star=2x_{1/2}-x$. Then $x^\star$ is the IMLA update
\eqref{eq:imla} started from $x$. The following lemma, whose elementary proof will not be given, shows that the
scheme is well defined under reasonable hypotheses.

\begin{lemma} [Well-posedness]
\label{lem:wellposed}

\begin{itemize}
\item
If $U\in \mathcal{C}^1(\mathbb{R}^d)$ and $U$ is convex, then for every $\delta>0$ \eqref{eq:midpoint} admits a unique solution.

\item If $U \in \mathcal{C}^{2}(\mathbb{R}^{d})$ and the eigenvalues of the Hessian of $U$ are bounded below by $-M$ for some $M>0$, then \eqref{eq:midpoint} admits a unique solution whenever $M\delta<2$.
\end{itemize}
\end{lemma}

Throughout the paper we write
\begin{equation}\label{eq:h}
h=\sqrt{2\delta},
\end{equation}
so that \eqref{eq:imla} reads $x^\star=x-(h^2/2)\nabla U(x_{1/2})+h\,\xi$.

\subsection{Hamiltonian interpretation and Metropolis correction} \label{sec:hmc}

Introduce the momentum $p=\xi$ and define
\begin{equation}\label{eq:mid-hmc}
p^\star=p-h\,\nabla U\!\Big(\frac{x^\star+x}{2}\Big),\qquad
x^\star=x+\frac{h}{2}\,(p+p^\star),
\end{equation}
which is  precisely one step of length $h$ of the implicit
midpoint rule applied to the Hamiltonian system
\[
\dot x=\frac{\partial H}{\partial p}=p,\qquad \dot p=-\frac{\partial H}{\partial x}=-\nabla U(x),\qquad
H(x,p)=\tfrac12\norm{p}^2+U(x).
\]
Eliminating $p^\star$ recovers \eqref{eq:imla},
thus IMLA is a single time step Hamiltonian Monte Carlo proposal with complete
Gaussian momentum refreshment at each time step. The implicit midpoint
map $(x,p)\mapsto(x^\star,p^\star)$ is well defined under the hypotheses of
Lemma~\ref{lem:wellposed}, time reversible, and symplectic, hence volume
preserving \cite{sanz2018numerical};
these are exactly the properties required by  HMC Metropolis
proposals \cite{BRSS18},
which accept $(x^\star,p^\star)$ with probability
$\min\{1,\exp(-(H(x^\star,p^\star)-H(x,p)))\}$.

The energy error $H(x^\star,p^\star)-H(x,p)$ may be expressed without reference to the momentum variable. In fact, from
\eqref{eq:mid-hmc},
\[
\tfrac12\norm{p^\star}^2-\tfrac12\norm{p}^2
=\tfrac12\,(p^\star-p)^{\top}(p^\star+p)
=-(x^\star-x)^{\top}\nabla U\!\Big(\frac{x^\star+x}{2}\Big),
\]
and therefore $H(x^\star,p^\star)-H(x,p)=\Delta(x,x^\star)$, where for
$x,y\in\mathbb{R}^d$ we define
\begin{equation*}
\Delta(x,y)=U(y)-U(x)-(y-x)^{\top}\nabla U\!\Big(\frac{x+y}{2}\Big).
\end{equation*}
The \emph{Metropolis adjusted implicit midpoint Langevin algorithm}
(MAIMLA) is the Metropolis--Hastings chain that, from the current
state $x$, proposes $x^\star$ by \eqref{eq:imla} and accepts it with
probability
\begin{equation}\label{eq:accept}
\alpha_\delta(x,x^\star)=\min\big\{1,\exp(-\Delta(x,x^\star))\big\}.
\end{equation}
By construction the chain is reversible
with respect to $\pi$.

Two features of \eqref{eq:accept} are worth stressing. First, although
the proposal is implicit, no Jacobian determinant appears: volume
preservation of the symplectic map  removes the factor
$\det(I_d+\delta\theta\nabla^2U)$ that burdens the Metropolised
$\theta$-schemes \cite{HSR21}. Second, when $U$ is
quadratic the midpoint identity $U(y)-U(x)=(y-x)^{\top}\nabla U((x+y)/2)$
holds exactly, so $\Delta\equiv0$ and every proposal is accepted: MAIMLA
 is exact for Gaussian targets. By implication, it may be expected to enjoy high acceptance rates in regions of the state space where the target is approximately Gaussian.

\subsection{Optimal scaling}\label{sec:optimal_scaling}

The Hamiltonian interpretation of Section~\ref{sec:hmc} allows us to
study the high-dimensional behaviour of MAIMLA within the framework of
\cite{OptimalScalingSymmetry26},
which applies to Metropolis algorithms whose proposal is one step of a
reversible, volume-preserving Hamiltonian integrator followed by a
Metropolis accept/reject step.

Now let $x=(x_1,\ldots,x_d)$, where each $x_i\in\mathbb{R}^m$, and consider the product target
\begin{equation*}
    \pi_d(x)
    \propto
    \prod_{i=1}^d
    \exp\left[-g\left(\frac{x_i}{\lambda_i}\right)\right],
    \qquad \lambda_i>0.
\end{equation*}
The choice $\lambda_i=1$ gives the i.i.d. product case, while $\lambda_i=i^{-\kappa}$, with $\kappa>0$, gives a product whose components have decreasing scales. For a single factor, introduce the Hamiltonian function
\[
    H(x,p)=\frac12\norm{p}^2+g(x),
    \qquad x,p\in\mathbb{R}^m,
\]
and let $\Psi_h(x,p)=(x^\star,p^\star)$ denote one step, with step size $h$, of an integrator for the associated Hamiltonian differential equations. Its energy error is
\[
    \Delta(x,p;h)=H(x^\star,p^\star)-H(x,p).
\]

The following result summarises Theorems~3.1 and 3.2 of \cite{OptimalScalingSymmetry26}, which, as $d\uparrow \infty$, yield formulas for the limiting acceptance probability  and the corresponding squared-jumping-distance in terms of the leading-order energy error. The equivariance hypothesis holds for any reasonable integrator.

\begin{theorem}[Optimal scaling of one-step Hamiltonian proposals]
\label{thm:general_optimal_scaling}

 Assume that:
\begin{itemize}
\item  $\Psi_h$ is equivariant under rescaling
the variables of the Hamiltonian system, reversible, and volume preserving.

\item There exist $\nu>0$ and functions $a(x,p)$ and $r(x,p;h)$ such that
\begin{equation}\label{eq:energy_error_expansion}
    \Delta(x,p;h)
    =h^\nu a(x,p)+h^\nu r(x,p;h),
    \qquad
    \lim_{h\to0}r(x,p;h)=0,
\end{equation}
and there exists a function $D(x,p)$, integrable with respect to the probability density  $\propto\exp(-H(x,p))$, such that
\begin{equation*}
    \sup_{0<h\leq1}\frac{|\Delta(x,p;h)|^3}{h^{3\nu}}
    \leq D(x,p).
\end{equation*}

\item
 The $\lambda_i$ are such that for constants $\gamma,K>0$,
\begin{equation*}
    \lim_{d\to\infty} d^{-\gamma}\sum_{i=1}^d\lambda_i^{-2\nu}= K,
    \qquad
    \lim_{d\to\infty} d^{-\gamma}\max_{1\leq i\leq d}\lambda_i^{-2\nu}
    =0.
\end{equation*}

\item The step length is chosen as
\begin{equation}\label{eq:hamiltonian_step_scaling}
    h_d=\ell d^{-\gamma/(2\nu)},
    \qquad \ell>0,
\end{equation}
\end{itemize}

Then, at stationarity, the expected acceptance probability converges, as $d\uparrow \infty$, to
\begin{equation}\label{eq:limiting_acceptance_general}
    A(\ell)
    =2\Phi\left(-\frac{\ell^\nu\sqrt{K\Sigma}}{2}\right),
\end{equation}
where $\Phi$ is the standard normal cumulative distribution function and
\begin{equation}\label{eq:optimal_sigma_general}
    \Sigma=\mathbb{E}\left[a(x,p)^2\right]
\end{equation}
($x$ has density  $\propto\exp(-g(x))$ and $p\sim\mathcal{N}(0,I_m)$).

If, in addition, the position increment satisfies $x^\star-x=hp+o(h)$ in the mean-square sense of \cite[Condition~3]{OptimalScalingSymmetry26}, then the squared jumping distance
\begin{equation}\label{eq:sjd_definition}
    \mathcal{SJD}_d
    =\frac1d\sum_{i=1}^d
      \mathbb{E}\left[
        \norm{x_i^\star-x_i}^2
        \mathbbm{1}_{\{\mathcal{U}\leq\alpha_d\}}
      \right],
    \qquad \mathcal{U}\sim\operatorname{Unif}(0,1),
\end{equation}
where $\alpha_d$ denotes the Metropolis acceptance probability for the $d$-component proposal,
satisfies
\begin{equation}\label{eq:sjd_limit_general}
    \lim_{d\to\infty}d^{\gamma/\nu}\mathcal{SJD}_d
    =m\ell^2A(\ell).
\end{equation}

\end{theorem}

We  apply Theorem~\ref{thm:general_optimal_scaling} to MAIMLA. The implicit midpoint rule is equivariant, volume-preserving, reversible and second order, and its leading energy error in \eqref{eq:energy_error_expansion} has $\nu=3$  for sufficiently smooth targets (see \cite{OptimalScalingSymmetry26} for details). Using \eqref{eq:h}, equation \eqref{eq:hamiltonian_step_scaling} becomes
\begin{equation*}
    \delta_d=\frac{\ell^2}{2}d^{-\gamma/3}.
\end{equation*}
For the i.i.d. product target, $\gamma=K=1$, whereas for $\lambda_i=i^{-\kappa}$ one has
\begin{equation}\label{eq:scaled_product_gamma}
    \gamma=6\kappa+1,
    \qquad
    K=\frac{1}{6\kappa+1}.
\end{equation}
Thus, in the latter case,
\[
    \delta_d=\frac{\ell^2}{2}d^{-(2\kappa+1/3)}.
\]

The method-dependent information is contained in the constant $\Sigma$ from \eqref{eq:optimal_sigma_general} and is given in the following corollary.
\begin{corollary}
  For the implicit midpoint integrator underlying MAIMLA, assuming sufficient smoothness and integrability of $g$, this constant is
\begin{equation}
\Sigma_{\mathrm{MAIMLA}}
=\mathbb{E}\left[
    \frac1{96}\norm{\nabla^3g}_F^2
    +\frac1{64}\norm{\nabla(\Delta g)}^2
\right].
\label{eq:sigma_maimla}
\end{equation}
Here the expectation is with respect to the single-component target with density $\propto\exp(-g)$.
\end{corollary}
We note that $\Sigma_{\mathrm{MAIMLA}}$ depends only on third derivatives of $g$, consistently with the fact that MAIMLA accepts every proposal for a Gaussian target. A well known argument \cite{OptimalScalingSymmetry26} shows that,  to attain optimal efficiency, one should try to ensure values of the acceptance rate $\approx 0.574$. This value does not take into account that, for implicit algorithms like MAIMLA,  increasing the step length may increase the computational cost of finding a proposal; when this cost is taken into consideration the optimal value will be smaller that $\approx 0.574$, \cite{OptimalScalingSymmetry26}. 

\subsubsection{Numerical illustration of the scaling limit}

We next illustrate the conclusions of Theorem~\ref{thm:general_optimal_scaling} numerically for the scaled-product target with $g(x)=x^4/4$ and $\lambda_i=i^{-\kappa}$, which contains the i.i.d.\ product as the special case
  $\kappa=0$. Throughout, we use a grid of step sizes $\delta=\tfrac{\ell^2}{2}d^{-(2\kappa+1/3)}$, with $\ell$ chosen such that the limiting acceptance probabilities form approximately an equidistant grid of 25 points.
  We measure efficiency by the expected
  squared jumping distance ($\mathcal{SJD}_d$) defined in \eqref{eq:sjd_definition}.
For each choice of the step size, the chain is initialised at stationarity and run for a single MAIMLA step, and the acceptance probability and $\mathcal{SJD}$ are estimated from $2\times 10^{3}$ independent realisations.

\begin{figure}[p!]
    \centering
    \begin{subfigure}{0.49\textwidth}
        \centering
         \includegraphics[height=0.219\textheight,keepaspectratio]{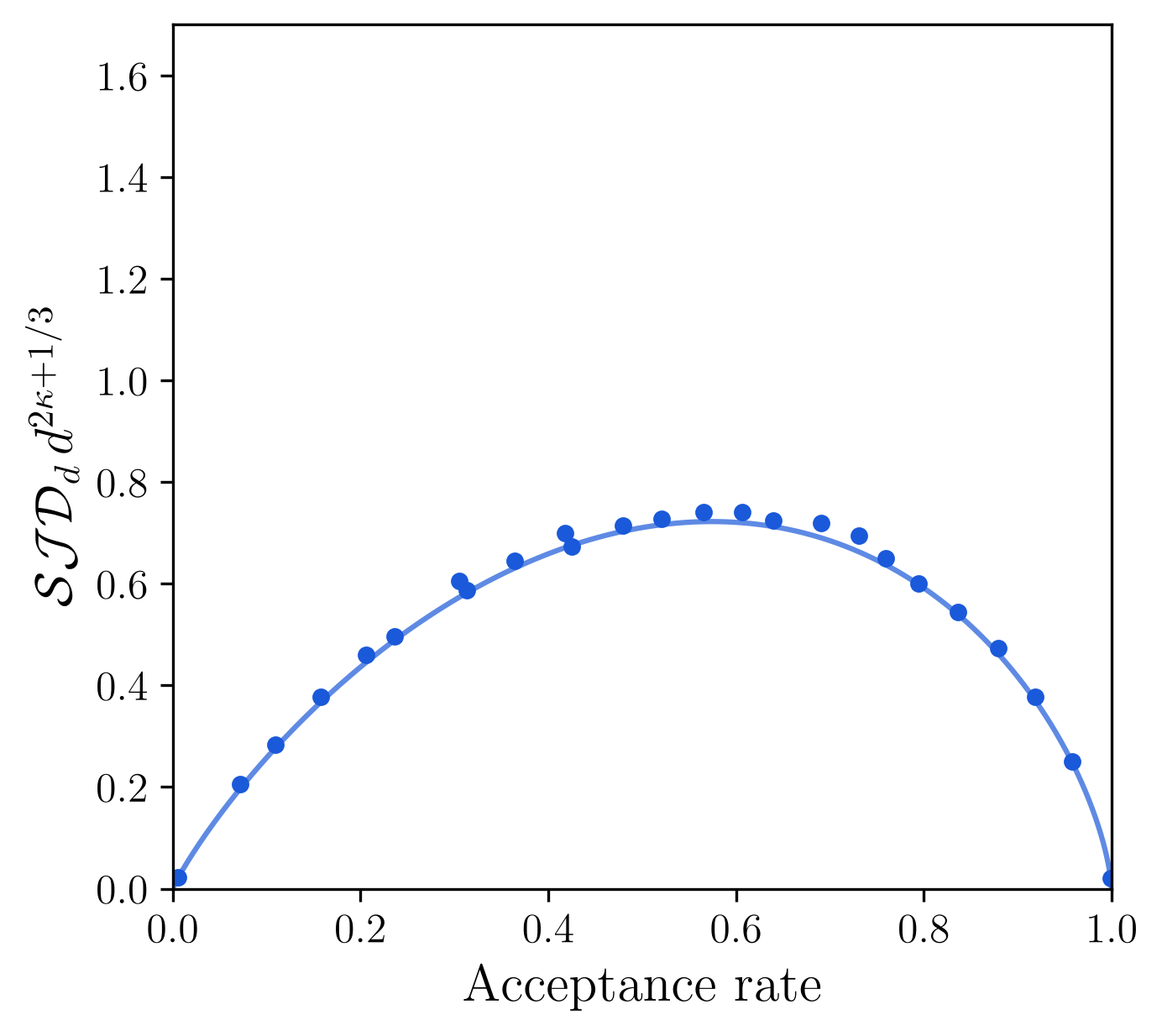}
        \includegraphics[height=0.218\textheight,keepaspectratio]{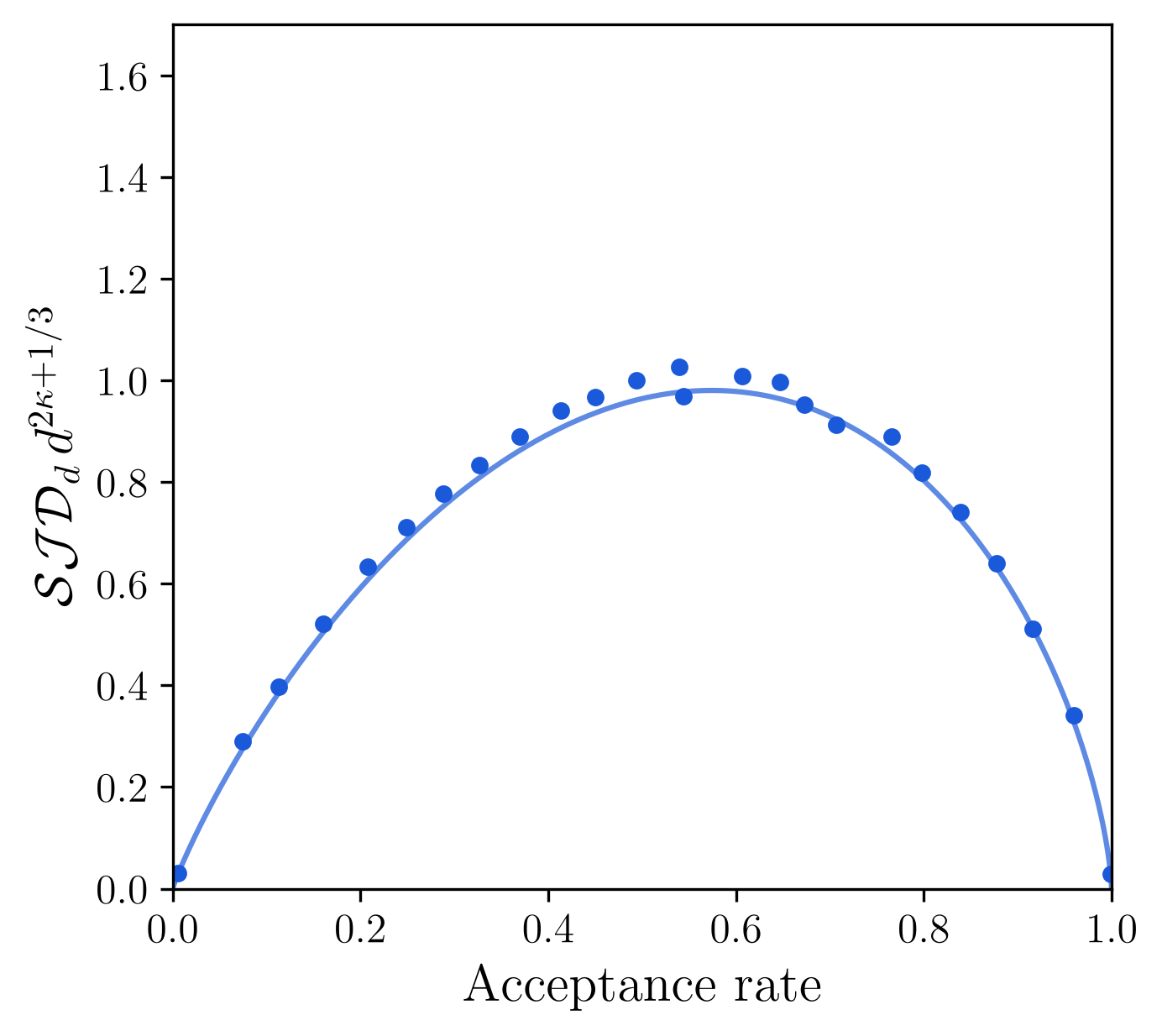}\\[0.25em]
        \includegraphics[height=0.218\textheight,keepaspectratio]{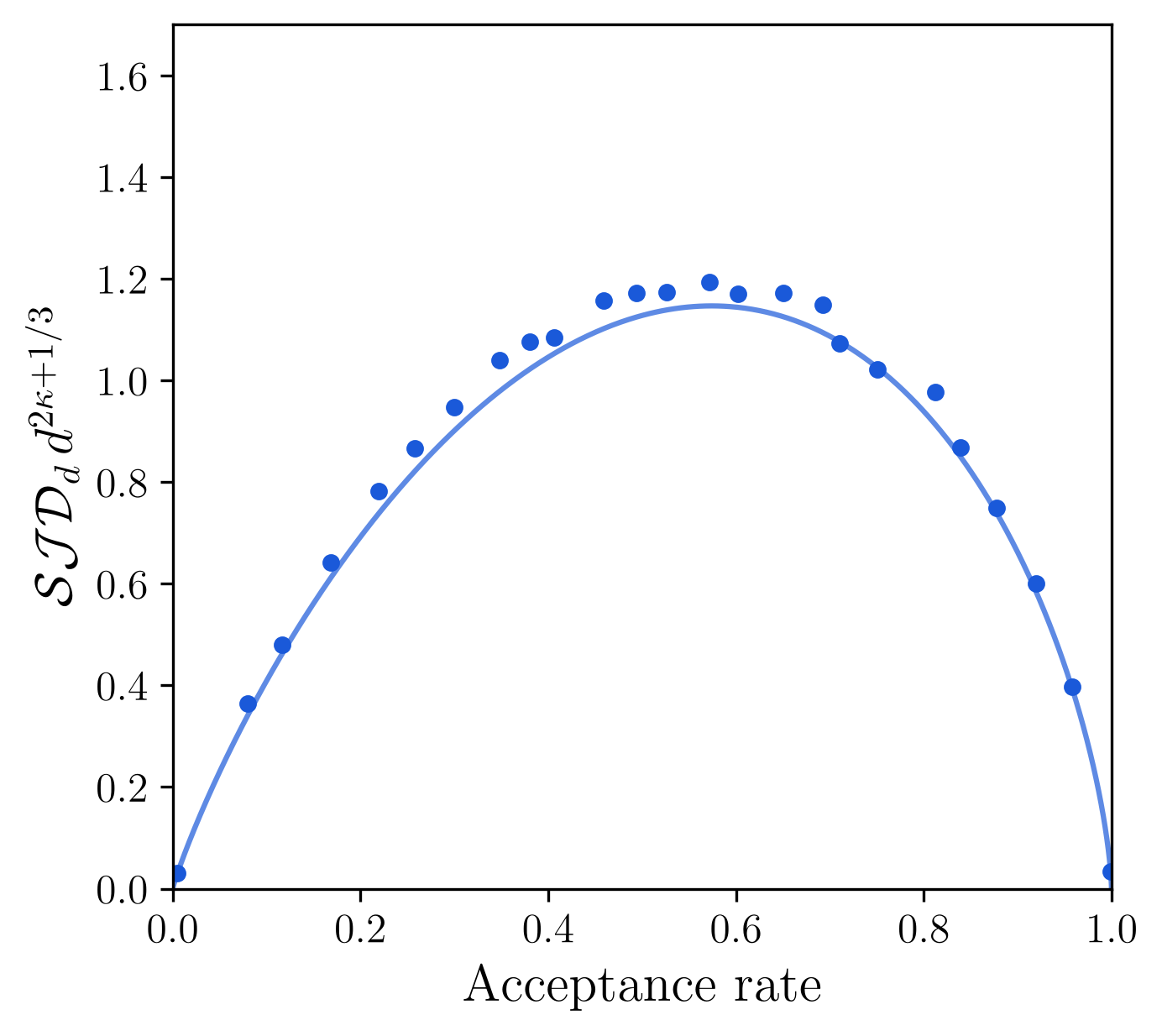}\\[0.25em]
        \includegraphics[height=0.218\textheight,keepaspectratio]{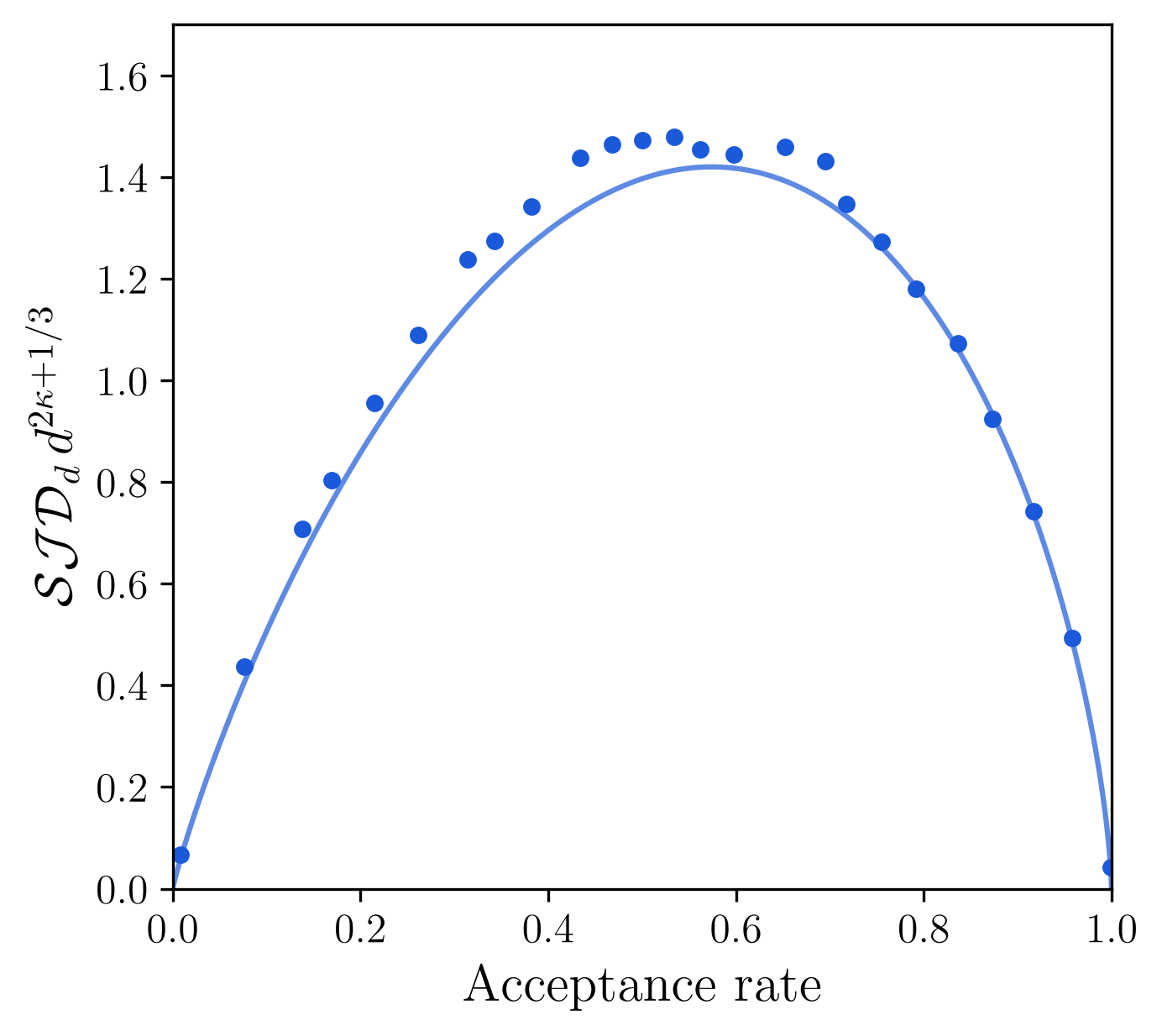}
        \subcaption{Scaled $\mathcal{SJD}_d$ against acceptance rate.}
        \label{fig:efficiency_accept_scaled}
    \end{subfigure}\hfill
    \begin{subfigure}{0.49\textwidth}
        \centering
        \includegraphics[height=0.219\textheight,keepaspectratio]{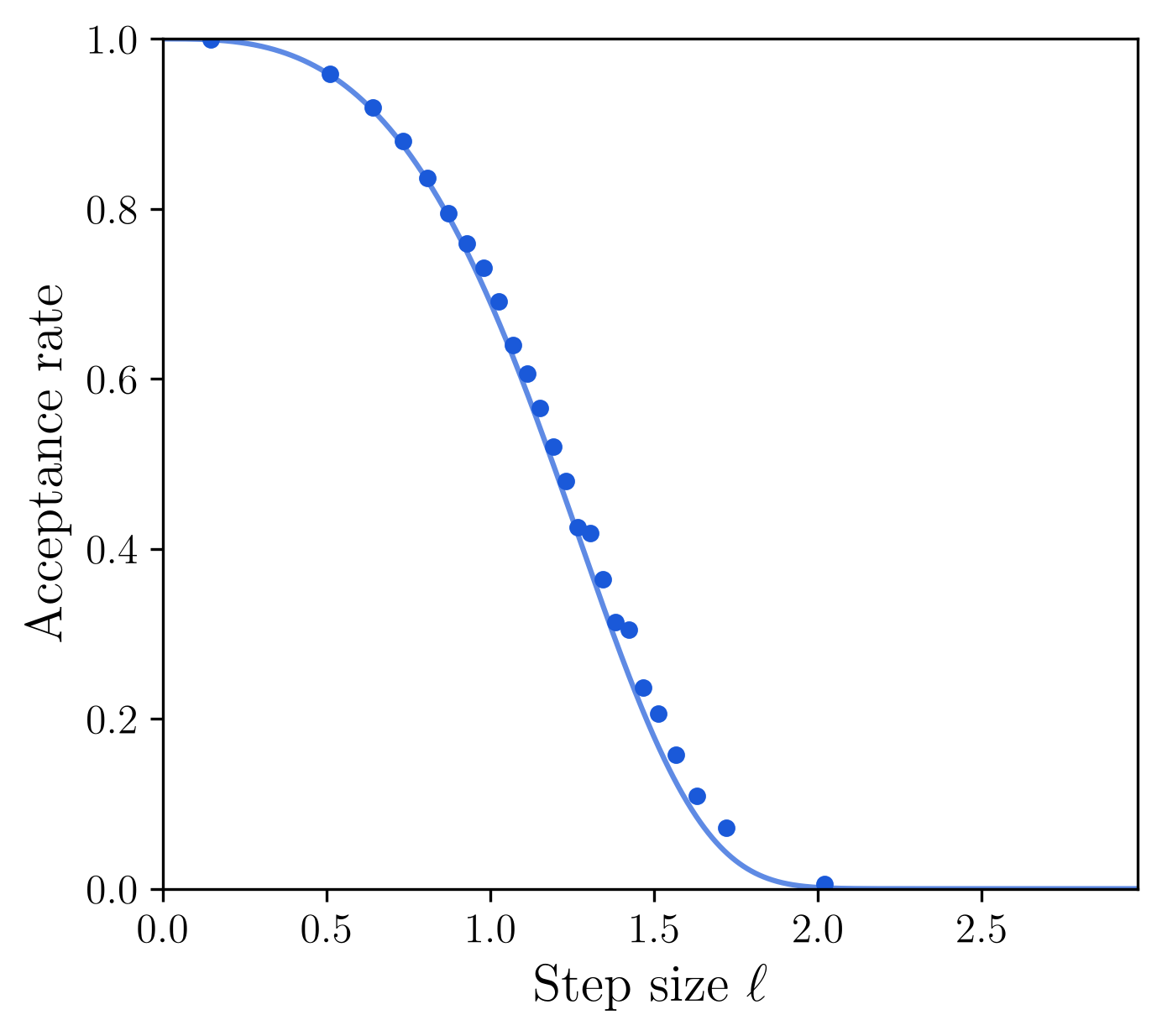}
        \includegraphics[height=0.218\textheight,keepaspectratio]{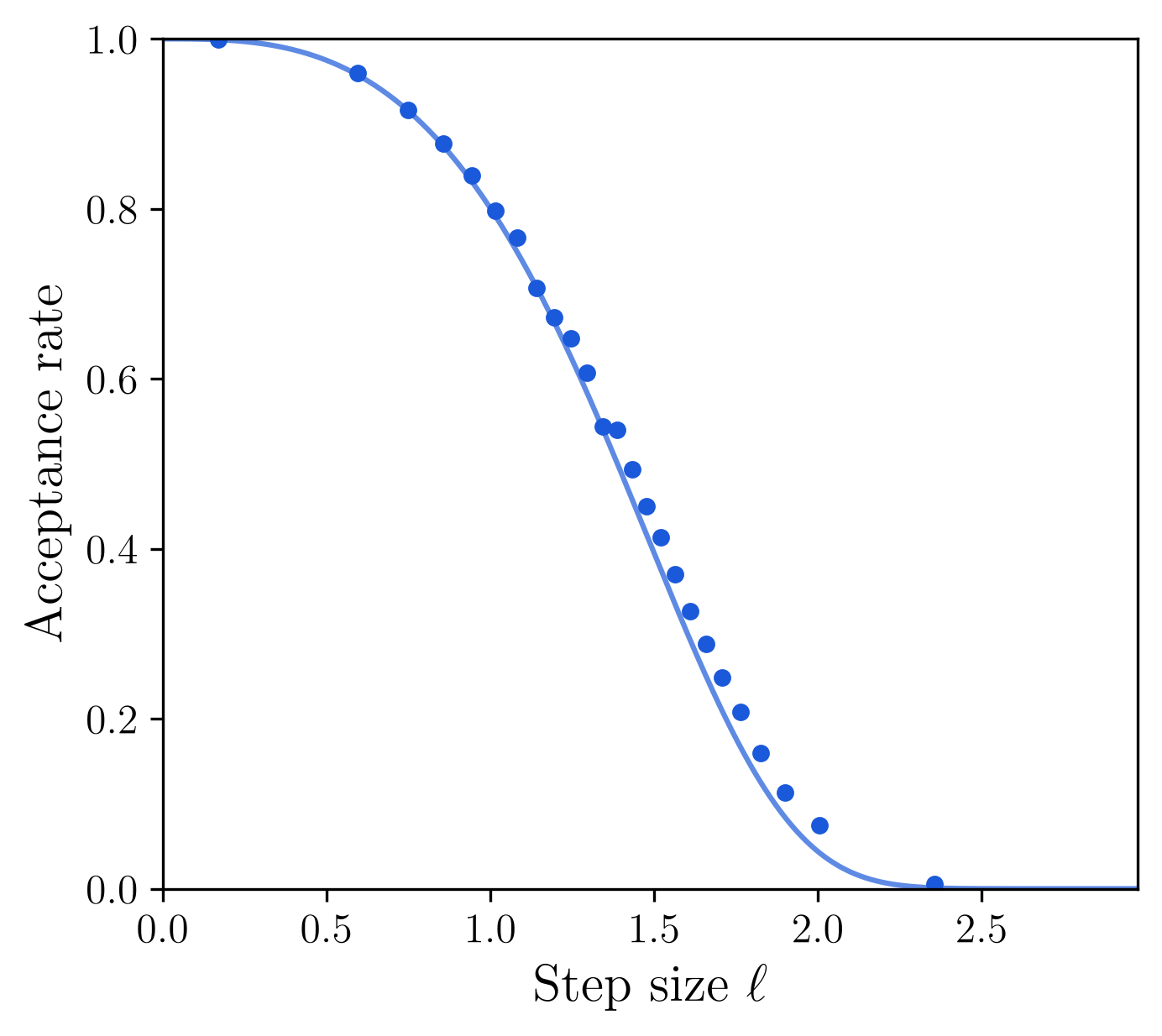}\\[0.25em]
        \includegraphics[height=0.218\textheight,keepaspectratio]{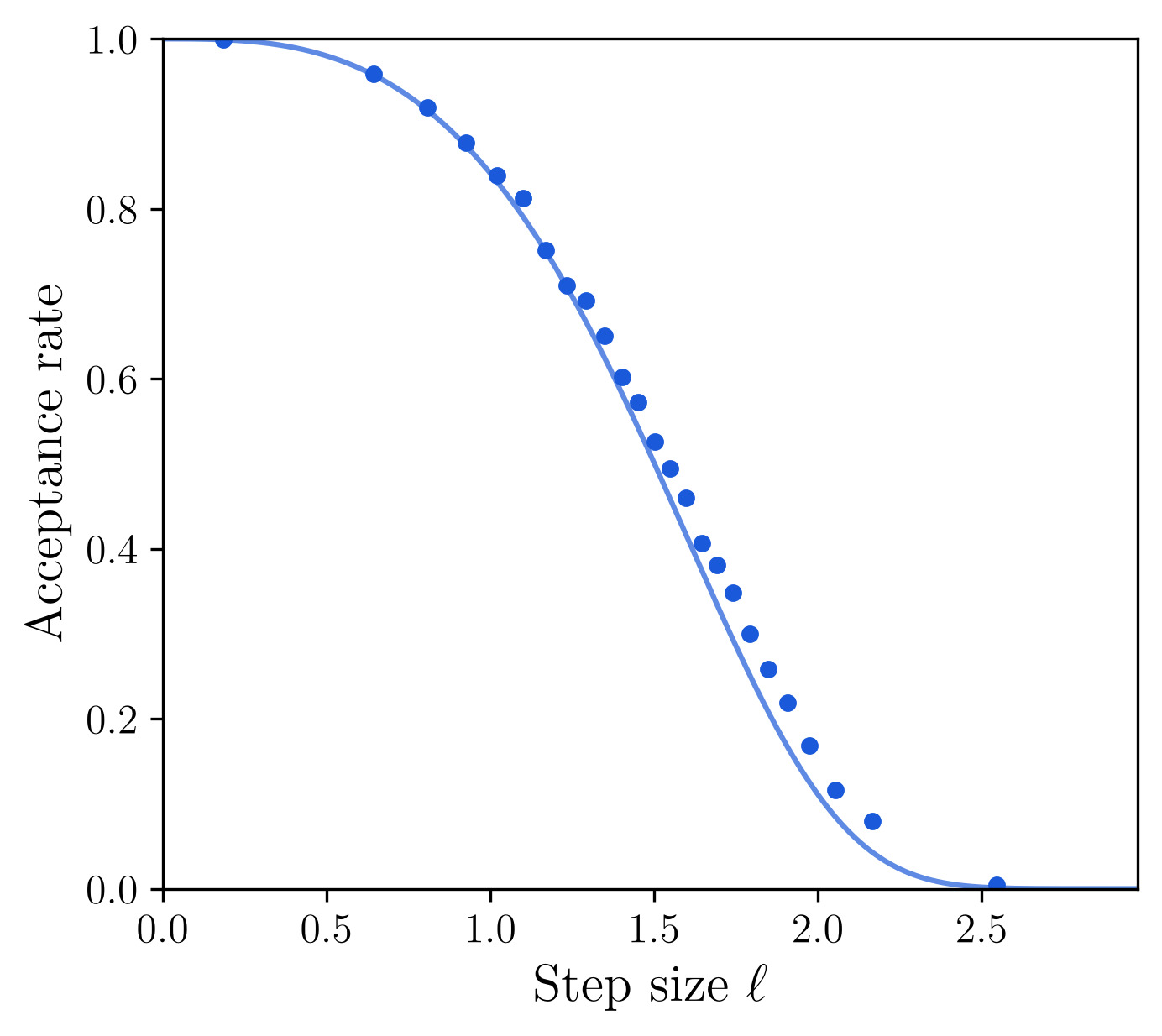}\\[0.25em]
        \includegraphics[height=0.218\textheight,keepaspectratio]{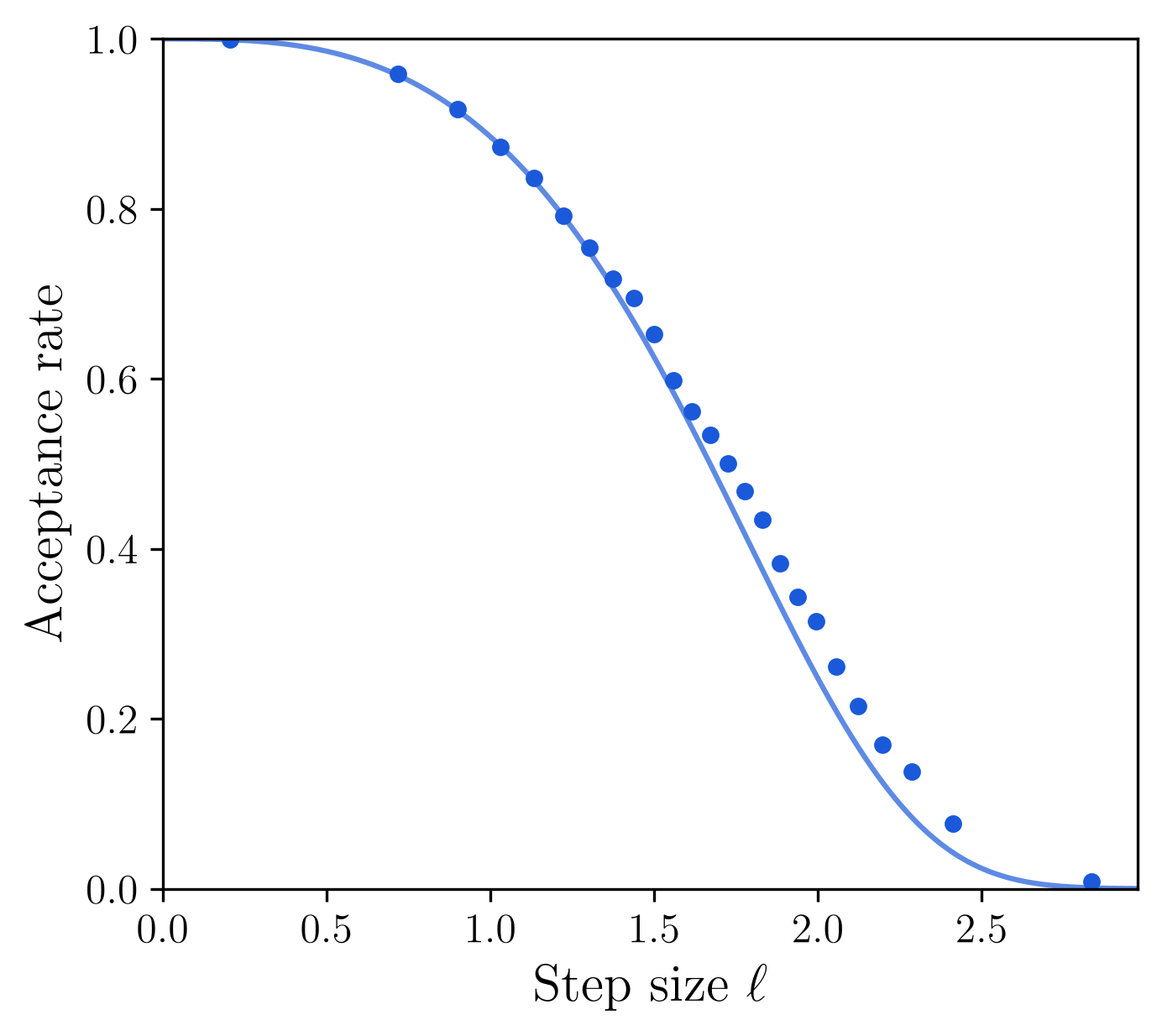}
        \subcaption{Acceptance rate against $\ell$.}
        \label{fig:ell_accept_scaled}
    \end{subfigure}
    \caption{Optimal-scaling behaviour of MAIMLA for the scaled quartic product target in dimension $d=10^{5}$, with $\kappa=0, 0.25,0.5,1.1$ from top to bottom.}
    \label{fig:IMLA_scaled_product}
\end{figure}

Figure~\ref{fig:IMLA_scaled_product} shows the i.i.d.\ case ($\kappa=0$) and the scaled product with $\kappa\in\{0.25,0.5,1.1\}$ at fixed $d=10^5$. The left column plots the $\mathcal{SJD}$ multiplied by $d^{2\kappa+1/3}$ against the empirical acceptance probability and the right column plots the acceptance probability against $\ell$, together with the limiting curves of Theorem~\ref{thm:general_optimal_scaling}. For  large $d $,
the finite-dimensional curves follow the asymptotic
predictions closely, confirming the exponent $\gamma=6\kappa+1$ in
\eqref{eq:scaled_product_gamma}.

Note that, in the left panels, the curves achieve its maximum when the acceptance rate is close to the value $0.574$ discussed above.

\subsubsection{Comparison with MALA}

The same Hamiltonian interpretation may be used for MALA, whose proposal corresponds to one step of the St\"ormer--Verlet integrator followed by a complete momentum refreshment. This integrator is also reversible and second order, so it has the same exponent $\nu=3$ and hence the same dimension-dependent step-size scaling as MAIMLA. Its method-dependent constant is
\begin{equation}
\Sigma_{\mathrm{MALA}}
=\mathbb{E}\left[
    \frac1{16}
    \norm{(\nabla^2g)\nabla g-\nabla(\Delta g)}^2
    +\frac1{24}\norm{\nabla^3g}_F^2
\right].
\label{eq:sigma_mala}
\end{equation}

We now  compare $\Sigma_{\mathrm{MAIMLA}}$ and $\Sigma_{\mathrm{MALA}}$
for one-dimensional factors $m=1$. In this case \eqref{eq:sigma_maimla} and \eqref{eq:sigma_mala} reduce to
\begin{equation}\label{eq:sigma_one_dimensional}
    \Sigma_{\mathrm{MALA}}
    =\frac1{48}\mathbb{E}\left[5g'''(x)^2+3g''(x)^3\right],
    \qquad
    \Sigma_{\mathrm{MAIMLA}}
    =\frac5{192}\mathbb{E}\left[g'''(x)^2\right].
\end{equation}
The following proposition summarises the comparison.
\begin{proposition}
\label{prop:sigma_comparison}
Let $m=1$ and let $g\in\mathcal{C}^3(\mathbb{R})$ then
\begin{equation}\label{eq:sigma_comparison_general}
    \Sigma_{\mathrm{MAIMLA}}\;\leq\;\frac58\,\Sigma_{\mathrm{MALA}}.
\end{equation}
If, in addition, $g$ is convex, then
\begin{equation}\label{eq:sigma_comparison_convex}
    \Sigma_{\mathrm{MAIMLA}}\;\leq\;\frac14\,\Sigma_{\mathrm{MALA}}.
\end{equation}
\end{proposition}
\begin{proof}
We first show that
\begin{equation}\label{eq:g2cubed_bound}
    \mathbb{E}\big[g''(x)^3\big]\;\geq\;-\,\mathbb{E}\big[g'''(x)^2\big].
\end{equation}
Integration by parts against the density $e^{-g}$ gives
$\mathbb{E}[f'(x)]=\mathbb{E}[f(x)g'(x)]$ for any $\mathcal{C}^1$ function $f$ with
$f e^{-g}\to0$ at $\pm\infty$. Taking $f=(g'')^2g'$, for which
$f'=2g'g''g'''+(g'')^3$, we obtain
\[
    \mathbb{E}\big[g''(x)^3\big]
    =\mathbb{E}\big[g'(x)^2g''(x)^2\big]-2\,\mathbb{E}\big[g'(x)g''(x)g'''(x)\big].
\]
Set $A=\mathbb{E}[g'(x)^2g''(x)^2]\geq0$ and $B=\mathbb{E}[g'''(x)^2]\geq0$. By the
Cauchy--Schwarz inequality, $|\mathbb{E}[g'g''\,g''']|\leq\sqrt{AB}$, and therefore
\[
    \mathbb{E}\big[g''(x)^3\big]\;\geq\;A-2\sqrt{AB}
    \;=\;\big(\sqrt A-\sqrt B\big)^2-B\;\geq\;-B,
\]
which is \eqref{eq:g2cubed_bound}.
Inserting \eqref{eq:g2cubed_bound} into the expression for $\Sigma_{\mathrm{MALA}}$
in \eqref{eq:sigma_one_dimensional} gives
\[
    \Sigma_{\mathrm{MALA}}
    =\frac1{48}\,\mathbb{E}\big[5g'''(x)^2+3g''(x)^3\big]
    \;\geq\;\frac{2}{48}\,\mathbb{E}\big[g'''(x)^2\big]
    =\frac{2}{48}\cdot\frac{192}{5}\,\Sigma_{\mathrm{MAIMLA}}
    =\frac85\,\Sigma_{\mathrm{MAIMLA}},
\]
which is \eqref{eq:sigma_comparison_general}. If, in addition, $g$ is convex, then $g''\geq0$, so that
$\mathbb{E}[g''(x)^3]\geq0$ and \eqref{eq:sigma_one_dimensional} yields
$\Sigma_{\mathrm{MALA}}\geq\frac{5}{48}\mathbb{E}[g'''(x)^2]=4\,\Sigma_{\mathrm{MAIMLA}}$,
which is \eqref{eq:sigma_comparison_convex}.
\end{proof}

In addition, it follows from \eqref{eq:limiting_acceptance_general} that, for the same step size, MAIMLA has a larger limiting acceptance probability. Alternatively, fixing a limiting acceptance probability $A$, equations \eqref{eq:limiting_acceptance_general} and \eqref{eq:sjd_limit_general} give
\begin{equation}\label{eq:sjd_as_function_of_acceptance}
    \lim_{d\to\infty}d^{\gamma/3}\mathcal{SJD}_d
    =\frac{2^{2/3}}{(K\Sigma)^{1/3}}
      A\left[\Phi^{-1}\left(1-\frac{A}{2}\right)\right]^{2/3}.
\end{equation}
Thus MAIMLA also has a larger limiting squared jumping distance when the two algorithms are tuned to the same acceptance probability. In such a regime  the limiting efficiency of MAIMLA exceeds that
of MALA by a factor of at least $(8/5)^{1/3}\approx1.17$ for every admissible
$g$, and by at least $4^{1/3}\approx1.59$ for convex $g$; the factor is
unbounded as $g$ approaches a quadratic. 

We illustrate this comparison using an i.i.d. product target whose one-dimensional factor is an equally weighted mixture of two Gaussian distributions with means $a$ and $-a$ and unit variance. Up to an additive constant, the corresponding potential is
\[
    g_a(x)=\frac12(x-a)^2-\log\left(1+e^{-2ax}\right).
\]
\begin{figure}[ht]
    \centering
    \begin{subfigure}[t]{0.49\textwidth}
        \includegraphics[width=0.9\linewidth]{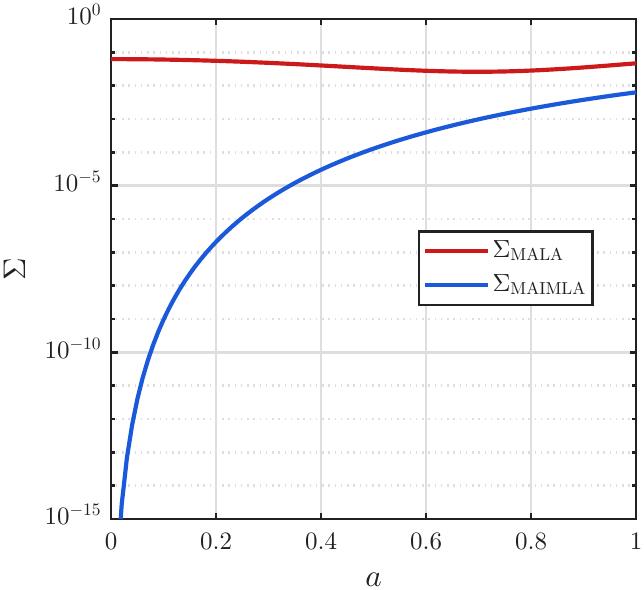}
        \subcaption{The constants $\Sigma_{\mathrm{MAIMLA}}$ and $\Sigma_{\mathrm{MALA}}$ as functions of $a$.}
        \label{fig:compa}
    \end{subfigure}
    \begin{subfigure}[t]{0.49\textwidth}
        \includegraphics[width=0.9\linewidth]{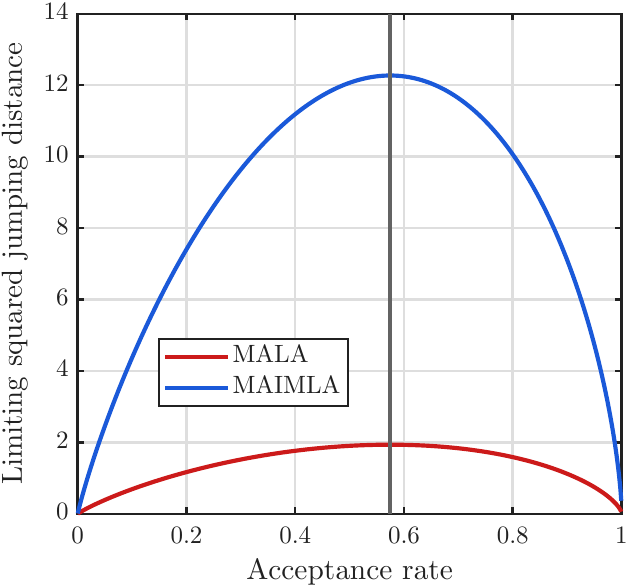}
        \subcaption{Limiting squared jumping distance against acceptance probability for $a=1/2$.}
        \label{fig:compb}
    \end{subfigure}
    \caption{Comparison of MAIMLA and MALA for the i.i.d. Gaussian-mixture product target with potential $g_a(x)=\frac12(x-a)^2-\log(1+e^{-2ax})$.}
    \label{fig:comparison_mala_maimla}
\end{figure}
We restrict $a\in[0,1]$, for which the potential is convex  and thus Lemma~\ref{lem:wellposed} applies. We evaluate $\Sigma_{\text{MAIMLA}}$ and $\Sigma_{\text{MALA}}$ using  \eqref{eq:sigma_one_dimensional} and plot in the left panel of Figure~\ref{fig:comparison_mala_maimla} their values as a function of $a$. As we can see $\Sigma_{\text{MAIMLA}}$ is always smaller than $\Sigma_{\text{MALA}}$ and as expected converges to zero when $a \rightarrow 0$, reflecting the exactness of the implicit midpoint rule for quadratic potentials. For the right panel, we set $a=1/2$ and compare the limiting square jumping distance for these two algorithms. As  predicted by \eqref{eq:sjd_as_function_of_acceptance} both algorithms achieve their fastest mixing for the same acceptance probability $0.574$, but the difference in the value of  $\Sigma$ leads to MAIMLA achieving values of  $\mathcal{SJD}_d$ six times as large as those provided by MALA.

Of course a fair comparison between MALA and MAIMLA must take into account the computational cost of generating the proposals, which, for implicit algorithms, depends heavily on the  target-dependent cost of solving the implicit equations; we shall return to this point later.

\section{Geometric ergodicity}
\label{sec:geo}
In this section we study the geometric ergodicity of IMLA and MAIMLA.
We begin by fixing terminology and recalling the two tools on which all
the proofs rest.

\paragraph{Geometric ergodicity and the Foster--Lyapunov criterion.}
Let $P$ be a Markov kernel. Recall that a set $C$ is \emph{small} for $P$ if there exist $\varepsilon>0$ and a
probability measure $\nu$ with $P(x,\cdot)\ge\varepsilon\,\nu(\cdot)$
for all $x\in C$.
We shall use the Foster--Lyapunov criterion in the
following form, see \cite[Theorem~15.0.1]{MT93}:
if $P$ is $\psi$-irreducible and aperiodic, every compact set is small,
$PV$ is finite everywhere, and
\begin{equation}\label{eq:drift}
\limsup_{\norm{x}\to\infty}\frac{PV(x)}{V(x)}<1 ,
\end{equation}
then $P$ has a unique invariant probability measure $\mu$, $\mu(V)<\infty$,
and the following \emph{$V$-geometric ergodicity} estimate holds  for suitable $K \in (0,\infty)$, $\rho\in(0,1)$:
\begin{equation*}%
\norm{P^n(x,\cdot)-\mu}_V\le K\,V(x)\,\rho^n ,\qquad x\in\mathbb{R}^d,\ n\ge0.
\end{equation*}

Condition \eqref{eq:drift} is equivalent
to the more familiar drift inequality $PV\le\lambda V+b\,\mathbbm{1}_C$
with $\lambda<1$ and $C$ compact. For the Metropolis--Hastings kernel of MAIMLA, whose proposal density and target $\pi$ are positive and continuous,
irreducibility and aperiodicity hold and every compact set is small by
\cite[Lemma~1.2]{MT96}
(see also \cite[Theorem~2.2]{RT96a}).
IMLA is not a Metropolis--Hastings chain, but the same conclusions
follow directly from the positivity and continuity of its transition
density: this is the content of Lemma~\ref{lem:minorisation_IMLA} below, whose
minorising measure is Lebesgue measure restricted to a ball and
therefore also yields Lebesgue-irreducibility and strong aperiodicity
\cite[Section~5.4]{MT93}.

\paragraph{Metropolis--Hastings kernels.}
Let $Q^{\text{IMLA}}_{\delta}$ denote the IMLA kernel, which under the conditions in Lemma~\ref{lem:wellposed} has
a density $q(x,\cdot)$, and let $P^{\text{MAIMLA}}_{\delta}$ denote the MAIMLA kernel. By
construction,
\begin{equation}\label{eq:MH-kernel}
P^{\text{MAIMLA}}_{\delta}(x,dy)=\alpha_\delta(x,y)\,Q^{\text{IMLA}}_{\delta}(x,dy)+\big(1-Q^{\text{IMLA}}_{\delta}\alpha_\delta(x)\big)\,\mathbbm{1}_x(dy),
\end{equation}
where $\mathbbm{1}_x$ is the Dirac mass at $x$ and $\alpha_\delta$ is
given by \eqref{eq:accept}; $1-Q^{\text{IMLA}}_{\delta}\alpha_\delta(x)$ is the probability of
rejecting a proposal from $x$. Two well-known facts about
\eqref{eq:MH-kernel} guide the analysis below. First, a Lyapunov
function for the unadjusted kernel $Q^{\text{IMLA}}_{\delta}$ is also a Lyapunov function for
$P^{\text{MAIMLA}}_{\delta}$ provided the rejection probability is bounded away from one outside
a compact set \cite[Theorem~3.1]{ST99};
second, and conversely, if the rejection probability is not bounded
away from one,
then $P^{\text{MAIMLA}}_{\delta}$ is not geometrically ergodic
\cite[Theorem~5.1]{RT96}.

\paragraph{Outline and heuristics.}
Geometric ergodicity of the unadjusted algorithm requires weaker
conditions than that of its Metropolised version, and we treat IMLA
first and MAIMLA second. The analysis of both rests on the
behaviour of the midpoint value $x_{1/2}$ of \eqref{eq:midpoint} when $\norm{x}$
is large. If $\norm{\nabla U(y)}\asymp\norm{y}^p$ with $p>1$, then for large
$\norm{x}$ the term $\frac{\delta}{2}\nabla U(x_{1/2})$ dominates the
left-hand side of \eqref{eq:midpoint} and
$\norm{x_{1/2}}\asymp\norm{x}^{1/p}\ll\norm{x}$; the precise statements are
Lemmas~\ref{lem:midpoint_bounds} and~\ref{lem:radial_inner_bound} in
Appendix~\ref{sec:midpoint}. Two consequences follow. For IMLA,
$\norm{X_1}^2-\norm{x}^2=-2\delta\,x_{1/2}^{\top}\nabla U(x_{1/2})+\cdots$
is strongly negative, which leads to an exponential Lyapunov function.
For MAIMLA, the proposal $x^\star=2x_{1/2}-x\approx-x$ reflects a tail
state through the bulk of the distribution; whether such a proposal is
accepted depends on $U(-x)$ being comparable to $U(x)$, which is why a
radial-symmetry condition on the leading-order part of $U$ enters in the
superlinear regime (Assumption~\ref{ass:ass3}). When $\nabla U$ is
globally Lipschitz the midpoint value $x_{1/2}$ is no longer small, the proposal moves by
$O(\delta\norm{x})$, and the acceptance probability is controlled instead
through the third derivative of $U$ (Assumption~\ref{ass:third_derivative_decay}). This is directly related to the fact that, contrary to MALA, no second derivative terms appear in the acceptance probability for MAIMLA. It is precisely the presence of those second derivative terms in MALA that allows one to prove geometric ergodicity  with no need to control the third derivative \cite{DM22}.

\subsection{IMLA}
Using a synchronous coupling argument it was shown in \cite{KDAPSZ24} that IMLA is contractive in the Wasserstein-$2$ distance assuming $U$ is gradient Lipschitz and strongly convex. Here we will work with   weaker assumptions, namely that $U$ satisfies Assumptions \ref{ass:ass1} and \ref{ass:ass2}, which allow for superlinear growth of $\nabla U$.  Our approach here is based on a Foster--Lyapunov argument  as in \cite{ST99}.
\begin{theorem}[Geometric ergodicity of IMLA]\label{thm:geometric_ergodicity_imla}
Suppose that $U\in\mathcal{C}^2(\mathbb{R}^d)$ is convex and satisfies
Assumptions~\ref{ass:ass1} and~\ref{ass:ass2}. Then, for every $\delta>0$,
the IMLA kernel $Q^{\text{IMLA}}_{\delta}$ has a unique invariant probability measure
$\pi_\delta$, and there exist $\rho\in(0,1)$ and $K<\infty$ such that
\begin{equation*}
\norm{(Q^{\text{IMLA}}_{\delta})^n(x,\cdot)-\pi_\delta}_V\le K\,V(x)\,\rho^n,\qquad x\in\mathbb{R}^d,\ n\ge0,
\end{equation*}
where $V(x)=\exp(a\norm{x}^2/8)$, with $a$  the constant in \eqref{eq:strong_ball}.
\end{theorem}
\begin{proof}
By Lemma~\ref{lem:Lyapunov_condition_unadjusted} below, $Q^{\text{IMLA}}_{\delta}V\le G(\norm{\cdot})V$ with
$G(r)\to0$, so \eqref{eq:drift} holds and $Q^{\text{IMLA}}_{\delta}V$ is finite everywhere. By
Lemma~\ref{lem:minorisation_IMLA} every compact set is small, and, as noted
after \eqref{eq:drift}, the minorising measure may be chosen as the Lebesgue measure on a
ball, so $Q^{\text{IMLA}}_{\delta}$ is Lebesgue-irreducible and strongly aperiodic. The claim follows
from the Foster--Lyapunov criterion.
\end{proof}

This subsection concludes with the two lemmas used to prove the preceding theorem.

\begin{lemma}\label{lem:Lyapunov_condition_unadjusted}
Suppose that $U\in \mathcal{C}^2(\mathbb{R}^d)$ satisfies \eqref{eq:strong_ball} and \eqref{eq:V_growth_cond}. Let
$        V(x) = \exp(a\norm{x}^2/8).
$
There exists a function $G:[0,\infty)\to [0,\infty)$ such that for all $x$,
\begin{equation}\label{eq:Lyapunov_IMLA}
\mathbb{E}[V(X_1)]\le G(\norm{x}) V(x), \quad \text{ and }\lim_{r\to\infty} G(r) =0,
\end{equation}
where $X_1$ is given by \eqref{eq:imla} with $X_0=x$.
\end{lemma}

\begin{proof}
    To simplify notation, denote by $X_{\frac{1}{2}}$ the midpoint $\frac{1}{2}(X_1+x)$. We start with the observation that
    \begin{align*}
        \norm{X_1}^2-\norm{x}^2 &= (X_1+x)^{\top}(X_1-x) \\
        &= (X_1+x)^{\top}\left(-\delta\nabla U(X_{\frac{1}{2}})+\sqrt{2\delta}\xi\right)\\
        &= -2\delta X_{\frac{1}{2}}^{\top} \nabla U(X_{\frac{1}{2}}) +2\sqrt{2\delta}X_{\frac{1}{2}}^{\top} \xi.
    \end{align*}
    Applying \eqref{eq:strong_ball}, we have
    \begin{align*}
        \norm{X_1}^2-\norm{x}^2
        \leq -2\delta a\norm{X_{\frac{1}{2}}}^2-2\delta b
        +2\sqrt{2\delta}X_{\frac{1}{2}}^{\top} \xi.
    \end{align*}
    For the noise term, Young's inequality gives
    \begin{equation*}
        2\sqrt{2\delta} X_{\frac{1}{2}}^{\top} \xi
        \leq \delta a \norm{X_{\frac{1}{2}}}^2 + \frac{2}{a} \norm{\xi}^2,
    \end{equation*}
    which leads to the estimate
    \begin{align}\label{eq:energy_estimate}
        \norm{X_1}^2-\norm{x}^2
        \leq -\delta a\norm{X_{\frac{1}{2}}}^2-2\delta b
        + \frac{2}{a} \norm{\xi}^2.
    \end{align}
    Applying \eqref{eq:midpoint_lower_bound} from Lemma~\ref{lem:midpoint_bounds}, we have
    \begin{align*}
        \norm{X_{\frac{1}{2}}} \geq g\left(\norm{x+\sqrt{\frac{\delta}{2}}\xi}\right), \quad
        \text{where } g(t) = \min\left\{\frac{1}{2}t,
        \left(\left(\frac{1}{\delta c_p}t-\frac{C_p}{c_p}\right)^+\right)^\frac{1}{p}\right\}.
    \end{align*}
    Therefore \eqref{eq:energy_estimate} becomes
    \begin{align*}
        \norm{X_1}^2-\norm{x}^2
        \leq -\delta a g\left(\norm{x+\sqrt{\frac{\delta}{2}}\xi}\right)^2
        -2\delta b + \frac{2}{a} \norm{\xi}^2.
    \end{align*}

    Since $g$ has sublinear growth for $p>1$, this suggests that we cannot use a quadratic Lyapunov function. Therefore, we work with the exponential Lyapunov function $V$. Exponentiating gives
    \begin{align*}
        \frac{V(X_1)}{V(x)}\leq  \exp\left(
        -\frac{a^2\delta}{8} g\left(\norm{x+\sqrt{\frac{\delta}{2}}\xi}\right)^2
        -\frac{a\delta b}{4} + \frac{1}{4}\norm{\xi}^2 \right).
    \end{align*}

    It remains to bound
    \begin{equation*}
        \mathbb{E}\left[\exp\left(
        -\frac{a^2\delta}{8} g\left(\norm{x+\sqrt{\frac{\delta}{2}}\xi}\right)^2
        -\frac{a\delta b}{4} + \frac{1}{4}\norm{\xi}^2 \right)\right].
    \end{equation*}
    To do this, split into two cases: (i) when $\xi$ is small, namely $A =\{\xi: \sqrt{\frac{\delta}{2}}\norm{\xi} \leq \frac{1}{2} \norm{x}\}$, and (ii) when $\xi$ is large, namely $A^c$.

    First, consider $\xi\in A$. Then
    \begin{equation*}
        \norm{x+\sqrt{\frac{\delta}{2}}\xi}
        \geq \norm{x}-\sqrt{\frac{\delta}{2}}\norm{\xi}
        \geq \frac{1}{2}\norm{x}.
    \end{equation*}
    Since $g$ is increasing, $g\left(\norm{x+\sqrt{\frac{\delta}{2}}\xi}\right)\geq g(\norm{x}/2)$, and hence
    \begin{equation*}
        \mathbb{E}\left[1_A\exp\left(
        -\frac{a^2\delta}{8} g\left(\norm{x+\sqrt{\frac{\delta}{2}}\xi}\right)^2
        -\frac{a\delta b}{4} + \frac{1}{4}\norm{\xi}^2 \right)\right]
        \leq e^{-\frac{a^2\delta}{8}g(\norm{x}/2)^2-\frac{a\delta b}{4}}
        \mathbb{E}\left[1_A\exp\left( \frac{1}{4}\norm{\xi}^2 \right)\right].
    \end{equation*}
    Since $\xi$ is a $d$-dimensional standard normal random variable, the moment generating function of a $\chi_{d}^2$ random variable gives
    \begin{equation*}
        \mathbb{E}\left[\exp\left( \frac{1}{4}\norm{\xi}^2 \right)\right]  = 2^{\frac{d}{2}},
    \end{equation*}
    and therefore
     \begin{equation*}
        \mathbb{E}\left[1_A\exp\left(
        -\frac{a^2\delta}{8} g\left(\norm{x+\sqrt{\frac{\delta}{2}}\xi}\right)^2
        -\frac{a\delta b}{4} + \frac{1}{4}\norm{\xi}^2 \right)\right]
        \leq 2^{\frac{d}{2}}e^{-\frac{a^2\delta}{8}g(\norm{x}/2)^2-\frac{a\delta b}{4}}.
    \end{equation*}

    Next, consider $\xi\in A^c$. Using $g\geq 0$, we have
    \begin{equation*}
        \mathbb{E}\left[1_{A^c}\exp\left(
        -\frac{a^2\delta}{8} g\left(\norm{x+\sqrt{\frac{\delta}{2}}\xi}\right)^2
        -\frac{a\delta b}{4} + \frac{1}{4}\norm{\xi}^2 \right)\right]
        \leq e^{-\frac{a\delta b}{4}}\mathbb{E}\left[1_{A^c}\exp\left(\frac{1}{4}\norm{\xi}^2 \right)\right].
    \end{equation*}

    Since $\norm{\xi}^2 \sim \chi_d^2$, for any $r>0$ we can write
    \begin{equation*}
        \mathbb{E}\left[1_{\{\norm{\xi}^2\geq r^2\}} e^{\frac{1}{4}\norm{\xi}^2}\right]
        = \frac{1}{2^{\frac{d}{2}}\Gamma(\frac{d}{2})}\int_{r^2}^\infty e^{\frac{1}{4}t} t^{\frac{d}{2}-1}e^{-\frac{t}{2}}\,dt
        = 2^{\frac{d}{2}} \mathbb{P}(T\geq r^2),
    \end{equation*}
    where $T\sim \mathrm{Gamma}(\frac{d}{2},\frac{1}{4})$. Markov's inequality gives
    \begin{equation*}
        \mathbb{E}\left[1_{\{\norm{\xi}^2\geq r^2\}} e^{\frac{1}{4}\norm{\xi}^2}\right]
        = 2^{\frac{d}{2}} \mathbb{P}(T\geq r^2)
        \leq 2^{\frac{d}{2}} e^{-\frac{1}{8}r^2}M_T\left(\frac{1}{8}\right)
        = 2^d e^{-\frac{1}{8}r^2}.
    \end{equation*}
    On $A^c$, we have $\norm{\xi}^2>\norm{x}^2/(2\delta)$. Applying the preceding bound with $r^2=\norm{x}^2/(2\delta)$ yields
    \begin{equation*}
        \mathbb{E}\left[1_{A^c}\exp\left(
        -\frac{a^2\delta}{8} g\left(\norm{x+\sqrt{\frac{\delta}{2}}\xi}\right)^2
        -\frac{a\delta b}{4} + \frac{1}{4}\norm{\xi}^2 \right)\right]
        \leq 2^d e^{-\frac{a\delta b}{4}}e^{-\frac{1}{16\delta}\norm{x}^2}.
    \end{equation*}
 Combining the two estimates gives
    \begin{equation*}
        \mathbb{E}[V(X_1)]
        \leq \left(2^{\frac{d}{2}}e^{-\frac{a^2\delta}{8}g(\norm{x}/2)^2-\frac{a\delta b}{4}}
        + 2^d e^{-\frac{a\delta b}{4}}e^{-\frac{1}{16\delta}\norm{x}^2}\right)V(x).
    \end{equation*}
\end{proof}

\begin{lemma}[Minorisation for IMLA]\label{lem:minorisation_IMLA}
Assume that $U\in\mathcal{C}^2(\mathbb{R}^d)$ is convex. Then for every $R>0$ the set $C=\{x:\ \norm{x}\le R\}$ is small.
\end{lemma}

\begin{proof}
By Lemma~\ref{lem:wellposed}, the map $x_{\frac{1}{2}}\mapsto x_{\frac{1}{2}}+\frac{\delta}{2}\nabla U(x_{\frac{1}{2}})$ is a $\mathcal{C}^1$ bijection of $\mathbb{R}^d$ with Jacobian $I_d+\frac{\delta}{2}\nabla^2U(x_{\frac{1}{2}})$, and $x^\star=2x_{\frac{1}{2}}-x$. A change of variables in \eqref{eq:midpoint} therefore shows that the IMLA kernel has transition density
\begin{equation*}
q(x,x^\star)=
(4\pi\delta)^{-d/2}
\det\left(I_d+\frac{\delta}{2}\nabla^{2}U\left(\frac{x+x^\star}{2}\right)\right)
\exp\left(
-\frac{\norm{x^\star-x+\delta\nabla U\left(\frac{x+x^\star}{2}\right)}^{2}}{4\delta}
\right).
\end{equation*}
Fix $R>0$. It suffices to find a lower bound on $q(x,x^\star)$ for all $x,x^\star\in C$. Since $U$ is convex,
\begin{equation*}
    \det\left(I_d+\frac{\delta}{2}\nabla^2 U\left(\frac{x+x^\star}{2}\right)\right) \geq 1.
\end{equation*}
Moreover, since $U\in\mathcal{C}^2(\mathbb{R}^d)$, the function
\begin{equation*}
    (x,x^\star)\mapsto\norm{x^\star-x+\delta\nabla U\left(\frac{x+x^\star}{2}\right)}^2
\end{equation*}
is continuous and hence bounded on $C\times C$. Consequently there exists $c_1>0$ such that $q(x,x^\star)\geq c_1$ for all $x,x^\star\in C$. Define
\[
\nu(A) := \frac{\mathrm{Leb}(A\cap C)}{\mathrm{Leb}(C)},
\qquad
\varepsilon := c_1\, \mathrm{Leb}(C).
\]
Then, for every measurable $A$ and every $x\in C$, $Q^{\text{IMLA}}_{\delta}(x,A)\ge\int_{A\cap C}q(x,x^\star)\,dx^\star\ge\varepsilon\,\nu(A)$,
which completes the proof.
\end{proof}

\subsection{MAIMLA}
We now study the geometric ergodicity of MAIMLA in two settings: when $\nabla U$ may grow superlinearly and when $\nabla U$ is globally Lipschitz. In both settings, we follow the proof strategy of Theorem~3.1 of \cite{ST99}. We first recall its key idea: a Lyapunov function for the unadjusted algorithm can also provide a drift condition for the adjusted algorithm when the mean acceptance probability is bounded away from zero outside a compact set. By Lemma~\ref{lem:Lyapunov_condition_unadjusted}, there is a nonnegative Lyapunov function $V:\mathbb{R}^d\to[0,\infty)$ and a function $G:[0,\infty)\to[0,\infty)$, with $G(r)\to0$ as $r\to\infty$, such that
\begin{equation}\label{eq:kernel_Lyap_imla}
    Q^{\text{IMLA}}_{\delta}V(x)\leq G(\norm{x})V(x),
\end{equation}
where $Q^{\text{IMLA}}_{\delta}V(x):=\int V(y)\,Q^{\text{IMLA}}_{\delta}(x,dy)$. By the definition of $P^{\text{MAIMLA}}_{\delta}$,
\begin{equation*}
    P^{\text{MAIMLA}}_{\delta}V(x) = Q^{\text{IMLA}}_{\delta}V(x)- \int_{\mathbb{R}^d}(1-\alpha_\delta(x,y)) V(y)\,Q^{\text{IMLA}}_{\delta}(x,dy)+(1-Q^{\text{IMLA}}_{\delta}\alpha_\delta(x))V(x);
\end{equation*}
the integral being subtracted is nonnegative and  thus, by \eqref{eq:kernel_Lyap_imla},
\begin{equation*}
    P^{\text{MAIMLA}}_{\delta}V(x) \leq \bigl(G(\norm{x})+1-Q^{\text{IMLA}}_{\delta}\alpha_\delta(x)\bigr)V(x).
\end{equation*}
Consequently, if
\begin{equation}\label{eq:rejection_prob_small}
    \limsup_{\norm{x}\to\infty}\bigl(1-Q^{\text{IMLA}}_{\delta}\alpha_\delta(x)\bigr) \leq \varepsilon
\end{equation}
for some $\varepsilon<1$, then
\begin{equation*}
    \limsup_{\norm{x}\to\infty}\frac{P^{\text{MAIMLA}}_{\delta}V(x)}{V(x)} \leq \varepsilon.
\end{equation*}
Thus, \eqref{eq:rejection_prob_small} yields a geometric drift condition for MAIMLA outside a compact set. This condition is weaker than requiring the rejection probability to converge to zero, as in \cite[Theorem~3.1]{ST99}. Indeed, the unadjusted drift condition \eqref{eq:Lyapunov_IMLA} also gives, for every $\lambda\in(0,1)$, a constant $C_\lambda<\infty$ such that
\begin{equation*}
    Q^{\text{IMLA}}_{\delta}V(x) \leq \lambda V(x)+C_\lambda.
\end{equation*}
Choosing $\lambda<1-\varepsilon$ then gives a geometric drift bound for the adjusted kernel, since $C_\lambda/V(x)\to0$ as $\norm{x}\to\infty$.

\subsubsection{Superlinear growth}

As with IMLA, we establish geometric ergodicity using a Foster--Lyapunov approach. When compared with IMLA we need to impose the additional Assumption \ref{ass:ass3} to be able to establish geometric ergodicity. This  implies that the leading order term in the potential $U$ is radially symmetric, and allows us to show that  the function $V$ from Lemma~\ref{lem:Lyapunov_condition_unadjusted} is also a Lyapunov function for MAIMLA. As  we will  show later, this extra condition cannot be dropped in general  when $\nabla U$ has superlinear growth.

\begin{theorem}\label{thm:geometric_ergodicity_MAIMLA}
Suppose that $U\in\mathcal{C}^2(\mathbb{R}^d)$ is convex and satisfies
Assumption~\ref{ass:ass3}. Then MAIMLA is $V$-geometrically ergodic for every
step size $\delta>0$, with $V$ as in Theorem~\ref{thm:geometric_ergodicity_imla}.
\end{theorem}

\begin{proof}
    Note that \eqref{eq:radial_condition} implies that the assumptions of Theorem~\ref{thm:geometric_ergodicity_imla} hold. Then by Lemma~\ref{lem:rejection_prob_zero} below we have
    \begin{equation}\label{eq:rejection_prob_convergence}
        \lim_{\norm{x} \to \infty} 1- Q^{\text{IMLA}}_{\delta}\alpha_\delta(x) =0
    \end{equation}
    and by \cite[Theorem~3.1]{ST99} MAIMLA is geometrically ergodic.
    \end{proof}

\begin{remark}[One-dimensional super-Gaussian targets]\label{rem:supergaussian}
A direct consequence of Theorem~\ref{thm:geometric_ergodicity_MAIMLA}
is that MAIMLA is geometrically ergodic, for every step size $\delta>0$,
when targeting the one-dimensional super-Gaussian family \eqref{eq:supergaussian}.
For these same targets the
Euler--Maruyama chain is transient and MALA is not geometrically ergodic
for any $\delta>0$ \cite[Theorems~3.2(b) and~4.2]{RT96}.
\end{remark}

The following lemma was the key in the proof of the theorem.

    \begin{lemma}\label{lem:rejection_prob_zero}
        Suppose that $U\in \mathcal{C}^2$ is convex and satisfies Assumption~\ref{ass:ass3}. Then \eqref{eq:rejection_prob_convergence} holds.
    \end{lemma}

    \begin{proof}
    Since the acceptance probability is given by \eqref{eq:accept}, if $x^\star\sim Q^{\text{IMLA}}_{\delta}(x,\cdot)$, then
    \begin{equation*}
        1-Q^{\text{IMLA}}_{\delta}\alpha_\delta(x)
        = \mathbb{E}\left[1-e^{-\Delta(x,x^\star)^+}\right]
        \leq 1-\mathbb{P}(\Delta\leq 0).
    \end{equation*}
Using \eqref{eq:midpoint}
    we have
    \begin{equation*}
        \Delta(x,x^\star)
        = U(x^\star)-U(x)
        + \delta \norm{\nabla U(x_{\frac{1}{2}})}^2
        - \sqrt{2\delta}\, \xi^{\top}\nabla U(x_{\frac{1}{2}}).
    \end{equation*}
    By Young's inequality
    \begin{equation*}
        \Delta(x,x^\star)
        \leq U(x^\star)-U(x)
        + 2\delta \norm{\nabla U(x_{\frac{1}{2}})}^2
        + \norm{\xi}^2.
    \end{equation*}

    Note that $x^\star=-x+2x_{\frac{1}{2}}$. As in the proof of Lemma~\ref{lem:Lyapunov_condition_unadjusted}, for large $x$, the midpoint $x_{\frac{1}{2}}$ is smaller than $x$. It is therefore convenient to work with
    \begin{equation*}
        \theta = \frac{x}{\norm{x}}, \quad
        \eta = \frac{x_{\frac{1}{2}}}{\norm{x}}.
    \end{equation*}
    Then
    \begin{equation*}
        \norm{x^\star}^2
        = \norm{-x+2x_{\frac{1}{2}}}^{2}
        = \norm{\norm{x}(-\theta+2\eta)}^2
        = \norm{x}^2\left(\norm{\theta}^2-4\theta^{\top}\eta+4\norm{\eta}^2\right).
    \end{equation*}
    By the binomial expansion, we have
    \begin{equation*}
        (1-\nu)^{\frac{1+p}{2}}
        = 1-\frac{1+p}{2}\nu + R_{\frac{1+p}{2}}(\nu),
        \quad |R_{\frac{1+p}{2}}(\nu)| \leq K_p |\nu|^2,
    \end{equation*}
    for some $K_p$ that depends only on $p$. Applying this expansion with $\nu = 4\theta^{\top}\eta-4\norm{\eta}^2$ (which will be justified later as small) gives
    \begin{equation*}
        \norm{-x+2x_{\frac{1}{2}}}^{1+p}
        = \norm{x}^{1+p}
        \left(1-2(p+1)(\theta^{\top}\eta-\norm{\eta}^2)
        + R_{\frac{1+p}{2}}(4\theta^{\top}\eta-4\norm{\eta}^2)\right).
    \end{equation*}
    Using \eqref{eq:radial_condition}, we can write
    \begin{align*}
        U(x^\star)-U(x) &= U(-x+2x_{\frac{1}{2}})-U(x) \\
        &= C\norm{-x+2x_{\frac{1}{2}}}^{1+p} - C\norm{x}^{p+1}+U_0(-x+2x_{\frac{1}{2}})-U_0(x)\\
        &= -2C(p+1)(\theta^{\top}\eta-\norm{\eta}^2)\norm{x}^{1+p}
        + C R_{\frac{1+p}{2}}(4\theta^{\top}\eta-4\norm{\eta}^2)\norm{x}^{1+p}
        + U_0(-x+2x_{\frac{1}{2}})-U_0(x)\\
        &= -2C(p+1)(x^{\top}x_{\frac{1}{2}}-\norm{x_{\frac{1}{2}}}^2)\norm{x}^{p-1}
        + C R_{\frac{1+p}{2}}\left(\frac{4x^{\top}x_{\frac{1}{2}}-4\norm{x_{\frac{1}{2}}}^2}{\norm{x}^2}\right)\norm{x}^{1+p}
        \\&+ U_0(-x+2x_{\frac{1}{2}})-U_0(x).
    \end{align*}

    Combining all the bounds we have
    \begin{align*}
        \Delta(x,x^\star) \leq -L_x+R_{1,x}+R_{2,x}+R_{3,x}+R_{4,x}+R_{5,x}
    \end{align*}
    where
    \begin{align*}
        L_x &= 2C(p+1)x^{\top}x_{\frac{1}{2}}\norm{x}^{p-1}\\
        R_{1,x} &= 2C(p+1)\norm{x_{\frac{1}{2}}}^2\norm{x}^{p-1} \\
        R_{2,x} &= C R_{\frac{1+p}{2}}\left(\frac{4x^{\top}x_{\frac{1}{2}}-4\norm{x_{\frac{1}{2}}}^2}{\norm{x}^2}\right)\norm{x}^{1+p}\\
        R_{3,x} &= U_0(-x+2x_{\frac{1}{2}})-U_0(x)\\
        R_{4,x} &= 2\delta \norm{\nabla U(x_{\frac{1}{2}})}^2\\
        R_{5,x} &= \norm{\xi}^2
    \end{align*}
    It remains to prove that, on the event $A_M=\{\norm{\xi}\leq M\}$ for some fixed $M>0$,
    \begin{equation*}
        \liminf_{\norm{x}\to\infty} \frac{L_x}{\norm{x}^{p+\frac{1}{p}} } >0
        \quad \text{and} \quad
        \lim_{\norm{x}\to\infty} \frac{R_{i,x}}{\norm{x}^{p+\frac{1}{p}} } =0
    \end{equation*}
    for all $i=1,\ldots,5$. Indeed, this implies that, for each fixed $M$, there exists $R_M>0$ such that $\Delta(x,x^\star)\leq0$ on $A_M$ whenever $\norm{x}\geq R_M$. Hence, for any $\varepsilon\in(0,1)$, we can choose $M$ such that
    \begin{equation*}
        \liminf_{\norm{x} \to\infty} \mathbb{P}(\Delta(x,x^\star) \leq 0) \geq \mathbb{P}(A_M)>1-\varepsilon.
    \end{equation*}
    As this holds for every $\varepsilon\in (0,1)$ we have that \eqref{eq:rejection_prob_convergence} holds.

    In what follows, we write $a\lesssim b$ if $a\leq Kb$ for a constant $K>0$ independent of $x$ and $\xi \in A_{M}$.
    \textit{Convergence of $L_{x}$.}
        By Lemma~\ref{lem:radial_inner_bound} we have
    \begin{align*}
        \frac{x^{\top}x_{\frac{1}{2}}}{\norm{x}^{1+1/p}}
        &\geq \frac{\delta}{2}C(1+p)\frac{\norm{x_{\frac{1}{2}}}^{p+1}}{\norm{x}^{1+1/p}}
        + \frac{\norm{x_{\frac{1}{2}}}^2}{\norm{x}^{1+1/p}}
        - \frac{\delta K_p\norm{x_{\frac{1}{2}}}^p}{\norm{x}^{1+1/p}}
        - \frac{K_p\left(1+\norm{\xi}^{\frac{p}{p-1}}\right)}{\norm{x}^{1+1/p}},
    \end{align*}
    for a constant $K_p>0$ depending only on $p$, $C_0$, and $c_0$.

    By Lemma~\ref{lem:midpoint_bounds},
    \begin{align*}
        \norm{x_{\frac{1}{2}}}^2
        &\leq
        \max\left\{\left(\frac{4}{\delta C_1}\right)^{\frac{1}{p-1}},
        \left(\frac{4}{\delta C_1}(\norm{z}-\tfrac{\delta}{2}c_1)\right)^{\frac{1}{p}}\right\}^2 \\
        &\lesssim
        \left(1+(\norm{z}-\tfrac{\delta}{2}c_1)^{\frac{2}{p}}\right).
    \end{align*}
    Since $(a-b)^{2/p}\lesssim a^{2/p}+1$ for $a\geq0$, we have
    \begin{equation*}
        (\norm{z}-\tfrac{\delta}{2}c_1)^{\frac{2}{p}}
        \lesssim 1+\norm{z}^{\frac{2}{p}}.
    \end{equation*}
    Using $z=x+\sqrt{\frac{\delta}{2}}\xi$, we obtain
    \begin{equation*}
        \norm{z}^{\frac{2}{p}}
        \lesssim \norm{x}^{\frac{2}{p}}+\norm{\xi}^{\frac{2}{p}}.
    \end{equation*}
    Since $\norm{\xi}^{2/p}\leq M^{2/p}$, it follows that restricted to the event $A_M$
    \begin{align}
        \norm{x_{\frac{1}{2}}}^2
        &\lesssim 1
        + \norm{x}^{2/p} .\label{eq:growth_square}%
    \end{align}
In particular, on the event $A_M$,
\begin{equation*}
|\nu| = \frac{\bigl|4x^{\top}x_{\frac{1}{2}}-4\norm{x_{\frac{1}{2}}}^2\bigr|}{\norm{x}^2}
\le \frac{4\norm{x_{\frac{1}{2}}}}{\norm{x}}+\frac{4\norm{x_{\frac{1}{2}}}^2}{\norm{x}^2}
\lesssim \norm{x}^{\frac{1}{p}-1}+\norm{x}^{\frac{2}{p}-2}\xrightarrow[\norm{x}\to\infty]{}0 ,
\end{equation*}
since $p>1$. Hence there exists $R'_M>0$ such that $|\nu|\le\tfrac12$ whenever $\norm{x}\ge R'_M$ and $\xi\in A_M$, which justifies the use of the binomial expansion above: for $|\nu|\le\tfrac12$ the remainder satisfies $|R_{\frac{1+p}{2}}(\nu)|\le K_p\nu^2$ with $K_p=\tfrac12\sup_{|t|\le 1/2}\bigl|\tfrac{d^2}{dt^2}(1-t)^{\frac{1+p}{2}}\bigr|<\infty$.
    By a similar argument, restricted to the event $A_M$,
    \begin{align*}
        \frac{\norm{x_{\frac{1}{2}}}^p}{\norm{x}^{1+1/p}}
        &\lesssim \frac{1}{\norm{x}^{1+1/p}}
        + \frac{\norm{x}}{\norm{x}^{1+1/p}} \\
        &= \norm{x}^{-(1+1/p)}+\norm{x}^{-1/p}.
    \end{align*}
    Therefore
    \begin{equation} \label{eq:convergence_square_term}
        \lim_{\norm{x}\to\infty}
       \frac{\norm{x_{\frac{1}{2}}}^p}{\norm{x}^{1+1/p}}
        = 0.
    \end{equation}

    Moreover, on the event $A_M$
    \begin{equation*}
        \lim_{\norm{x}\to\infty}\frac{K_p\left(1+\norm{\xi}^{\frac{p}{p-1}}\right)}{\norm{x}^{1+1/p}} =0.
    \end{equation*}

    Therefore,
        \begin{align*}
        \liminf_{\norm{x}\to\infty}\frac{x^{\top}x_{\frac{1}{2}}}{\norm{x}^{1+1/p}}
        &\geq \liminf_{\norm{x}\to\infty}\frac{\delta}{2}C(1+p)\frac{\norm{x_{\frac{1}{2}}}^{p+1}}{\norm{x}^{1+1/p}}.
    \end{align*}
    Note that \eqref{eq:radial_condition} implies that both conditions of Lemma~\ref{lem:midpoint_bounds} are satisfied and therefore
    \begin{equation*}
        \norm{x_{\frac{1}{2}}}^{p+1} \geq   \min\left\{2^{-(p+1)}\norm{z}^{p+1},
            \left(\left(\frac{1}{\delta c_p}\norm{z}-\frac{C_p}{c_p}\right)^+\right)^{1+\frac{1}{p}}\right\}
    \end{equation*}
On the event $A_M$ and for $\norm{x}$ sufficiently large, the second term in the minimum is positive and active. Thus
 \begin{equation*}
        \norm{x_{\frac{1}{2}}}^{p+1}\geq
            \left(\frac{1}{\delta c_p}\norm{z}-\frac{C_p}{c_p}\right)^{1+\frac{1}{p}} , \quad z=x+\sqrt{\frac{\delta}{2}}\xi.
    \end{equation*}
Now $\norm{z} \geq \norm{x} - \sqrt{\frac{\delta}{2}} \norm{\xi} \geq \norm{x} - M\sqrt{\frac{\delta}{2}}$ for $\xi\in A_M$, hence
    \begin{equation*}
         \norm{x_{\frac{1}{2}}}^{p+1} \geq
            \left(\frac{1}{\delta c_p}\norm{x} -\frac{1}{\sqrt{2\delta}c_p}M -\frac{C_p}{c_p}\right)^{1+\frac{1}{p}} .
    \end{equation*}
    Therefore on the event $A_M$ we have
    \begin{equation*}
        \liminf_{\norm{x} \to\infty} \frac{\norm{x_{\frac{1}{2}}}^{p+1}}{\norm{x}^{1+1/p}} \geq \frac{1}{(\delta c_p)^{1+\frac{1}{p}}}>0,
    \end{equation*}
    which implies that
    \begin{equation*}
        \liminf_{\norm{x}\to\infty}\frac{L_x}{\norm{x}^{p+\frac{1}{p}} } >0.
    \end{equation*}

    \textit{Convergence of $R_{1,x}$.} This follows immediately from \eqref{eq:convergence_square_term}.

    \textit{Convergence of $R_{2,x}$.} Using that $R_{\frac{1+p}{2}}$ grows at most quadratically, we have
    \begin{equation*}
        \frac{R_{2,x}}{\norm{x}^{p+\frac{1}{p}} } \lesssim \frac{\left(4x^{\top}x_{\frac{1}{2}}-4\norm{x_{\frac{1}{2}}}^2\right)^2}{\norm{x}^{\frac{1}{p}+3} } \lesssim  \frac{\norm{x}^2\norm{x_{\frac{1}{2}}}^2+ \norm{x_{\frac{1}{2}}}^4}{\norm{x}^{\frac{1}{p}+3} }
    \end{equation*}
    By \eqref{eq:growth_square} we have on the event $A_M$
    \begin{equation*}
        \frac{R_{2,x}}{\norm{x}^{p+\frac{1}{p}} }  \lesssim  \frac{\norm{x}^2+\norm{x}^{2+2/p}+ 1+\norm{x}^{4/p}}{\norm{x}^{\frac{1}{p}+3} }.
    \end{equation*}
   The right-hand side converges to zero as $\norm{x}$ tends to infinity since $p>1$.

   \textit{Convergence of $R_{3,x}$.} By the assumptions on $U_0$, we have that $U_0(x)$ grows at most of order $\norm{x}^p$. Therefore
   \begin{equation*}
    \frac{|R_{3,x}|}{\norm{x}^{p+\frac{1}{p}}} \lesssim \frac{1+\norm{x}^p +\norm{x_{\frac{1}{2}}}^p}{\norm{x}^{p+\frac{1}{p}}}
   \end{equation*}
   By \eqref{eq:growth_square} on the event $A_M$ we have
   \begin{equation*}
    \lim_{\norm{x}\to\infty}\frac{|R_{3,x}|}{\norm{x}^{p+\frac{1}{p}}} \lesssim\lim_{\norm{x}\to\infty} \frac{1+\norm{x}^p +\norm{x}}{\norm{x}^{p+\frac{1}{p}}} =0.
   \end{equation*}

   \textit{Convergence of $R_{4,x}$.} Note that our assumptions imply that \eqref{eq:V_growth_cond} holds, and hence
   \begin{equation*}
    \frac{R_{4,x}}{\norm{x}^{p+\frac{1}{p}} }  \lesssim  \frac{\norm{\nabla U(x_{\frac{1}{2}})}^2}{\norm{x}^{p+\frac{1}{p}} } \lesssim   \frac{1+\norm{x_{\frac{1}{2}}}^{2p}}{\norm{x}^{p+\frac{1}{p}} }.
   \end{equation*}
   Applying \eqref{eq:growth_square} when restricted to the event $A_M$ we have
   \begin{equation*}
     \frac{R_{4,x}}{\norm{x}^{p+\frac{1}{p}} }  \lesssim   \frac{1+\norm{x}^{2}}{\norm{x}^{p+\frac{1}{p}} }.
   \end{equation*}
   Since $p>1$, the right-hand side converges to zero as $\norm{x}\to \infty$.

    \textit{Convergence of $R_{5,x}$.} This follows immediately on the event $A_M$.

\end{proof}

\begin{figure}[h]
    \centering
    \includegraphics[width=0.6\textwidth]{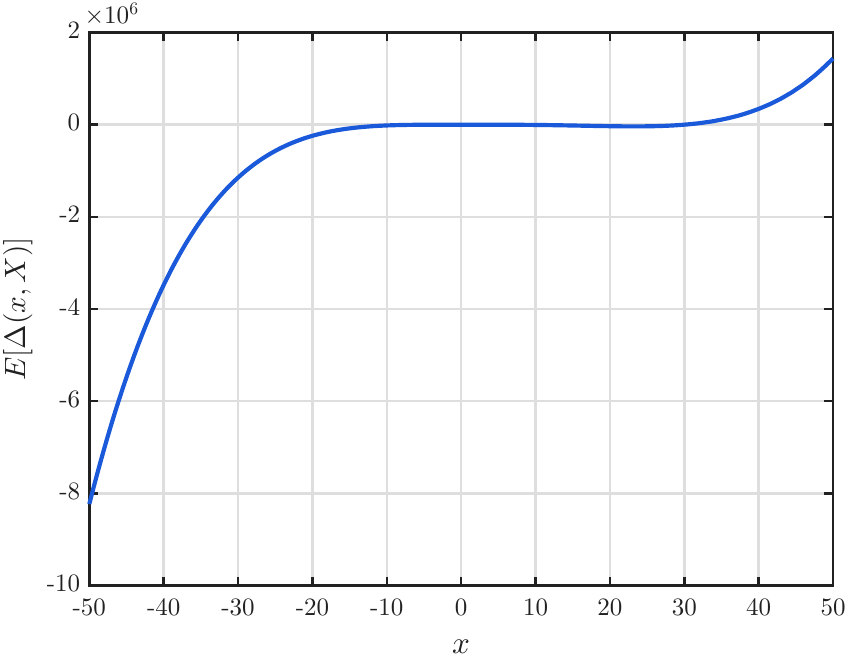}
    \caption{Expected value of $\Delta(x,x^\star)$ for the asymmetric one-dimensional potential obtained by averaging $100$ independent copies. Here $\delta=1$ and the implicit step was solved using Newton's method.}
    \label{fig:asymmetric_energy_estimate}
\end{figure}

\begin{remark} \label{rem:asym}
Note that, to establish geometric ergodicity for MAIMLA, we had to impose more restrictive assumptions than for IMLA. In particular, we assumed that the leading-order term of $U$ is radially symmetric. This symmetry assumption is vital for the proof, since for large $x$ the proposal is approximately $-x$. Therefore, if $U(x)$ differs greatly from $U(-x)$, the energy error $\Delta(x,x^\star)$ can be very large, and one can find values of $x$ for which the acceptance rate is arbitrarily close to zero. To demonstrate this, we consider a one-dimensional example where
\[
U(x) =
\begin{cases}
2x^4, & x < -1, \\[6pt]
\bigl(\frac{3}{2}- \tfrac{1}{2}\sin \left( \frac{\pi}{2} x\right)\bigr)x^4, & -1 \le x \le 1, \\[6pt]
x^4, & x > 1.
\end{cases}
\]
In Figure~\ref{fig:asymmetric_energy_estimate}, we plot $\mathbb{E}[\Delta(x,x^\star)]$ for a range of values of $x$. We observe that, although $\mathbb{E}[\Delta(x,x^\star)]$ is negative for negative $x$, it becomes extremely large when $x$ is large and positive, leading to a very low acceptance probability.
\end{remark}

\subsubsection{Lipschitz gradient}

We have observed that if $U$ is quadratic then $\Delta=0$ and hence MAIMLA accepts all proposals with probability one. As a consequence, if one considers a Taylor expansion of $\Delta$ for small $\delta$ we see that the terms involving the first and second order derivative of $U$ vanish. This implies that we will need to impose some assumption on the higher order derivatives. We will impose Assumption~\ref{ass:third_derivative_decay} which implies that the tails of the distribution $\pi$  approach those of a Gaussian distribution.

\begin{theorem}\label{thm:maimla-lip}
Suppose that $U\in\mathcal{C}^3(\mathbb{R}^d)$ is convex and satisfies Assumptions~\ref{ass:ass1}
and~\ref{ass:third_derivative_decay}, and let $0<\delta<2/L$. Then MAIMLA is
$V$-geometrically ergodic with $V(x)=\exp(a\norm{x}^2/8)$.
\end{theorem}
\begin{proof}
Assumption~\ref{ass:third_derivative_decay} gives Assumption~\ref{ass:ass2} with
$p=1$, so Lemma~\ref{lem:Lyapunov_condition_unadjusted} applies;
Lemma~\ref{lem:lip_acceptance} gives \eqref{eq:rejection_prob_small}; we can thus conclude
as in Theorem~\ref{thm:geometric_ergodicity_MAIMLA}.
\end{proof}
\begin{lemma}\label{lem:lip_acceptance}
Under the hypotheses of Theorem~\ref{thm:maimla-lip} there is $c_\delta>0$ with
$\inf_{x}Q^{\text{IMLA}}_{\delta}\alpha_\delta(x)\ge c_\delta$; in particular
\eqref{eq:rejection_prob_small} holds with $\varepsilon=1-c_\delta$.
\end{lemma}

\begin{proof}
 By Assumption~\ref{ass:ass1}, $U$ is coercive and hence has a minimiser, therefore to simplify the exposition, we assume that $\nabla U(0)=0$ which may be imposed by translating a minimiser to the origin.
Using a Taylor expansion, we obtain
\begin{equation*}
    \Delta(x,x^\star)= \int_{-\frac{1}{2}}^{\frac{1}{2}} \frac{1}{2}(\frac{1}{2}-|s|)^2 \nabla^3 U(x_{\frac{1}{2}}+s(x^\star-x))\left[x^\star-x,x^\star-x,x^\star-x\right] ds.
\end{equation*}
Using the proposal equation~\eqref{eq:imla}, we obtain
\begin{equation*}
    \Delta(x,x^\star) = -\delta^3 A_0 + 3\sqrt{2}\delta^{\frac{5}{2}} A_1-6\delta^2 A_2 + 2^{\frac{3}{2}}\delta^{\frac{3}{2}}A_3
\end{equation*}
where $y_s=x_{\frac{1}{2}}+s(x^\star-x)$ and
\begin{equation*}
    w(s)=\frac{1}{2}\left(\frac{1}{2}-|s|\right)^2.
\end{equation*}
Then
\begin{align*}
    A_0 &= \int_{-\frac{1}{2}}^{\frac{1}{2}} w(s)
    \nabla^3 U(y_s)\left[\nabla U(x_{\frac{1}{2}}),\nabla U(x_{\frac{1}{2}}),\nabla U(x_{\frac{1}{2}})\right] \, ds,\\
    A_1 &= \int_{-\frac{1}{2}}^{\frac{1}{2}} w(s)
    \nabla^3 U(y_s)\left[\nabla U(x_{\frac{1}{2}}),\nabla U(x_{\frac{1}{2}}),\xi\right] \, ds,\\
    A_2 &= \int_{-\frac{1}{2}}^{\frac{1}{2}} w(s)
    \nabla^3 U(y_s)\left[\nabla U(x_{\frac{1}{2}}),\xi,\xi\right] \, ds,\\
    A_3 &= \int_{-\frac{1}{2}}^{\frac{1}{2}} w(s)
    \nabla^3 U(y_s)\left[\xi,\xi,\xi\right] \, ds.
\end{align*}

Let us first consider $A_0$. Since $\nabla U$ is $L$-Lipschitz and $\nabla U(0)=0$, we have
\begin{equation}\label{eq:first_bound}
   \norm{\nabla U(x)} \leq L \norm{x}.
\end{equation}
Using~\eqref{eq:third_decay} and~\eqref{eq:first_bound} gives
\begin{align*}
    |A_0| \leq CL^3\int_{-\frac{1}{2}}^{\frac{1}{2}} \frac{1}{2}\left(\frac{1}{2}-|s|\right)^2
    \frac{\norm{x_{\frac{1}{2}}}^3}{(1+\norm{y_s})^3} \, ds.
\end{align*}

Note that
\begin{align*}
    \norm{y_s} \geq \norm{x_{\frac{1}{2}}}-|s|\norm{x^\star-x},
\end{align*}
and, by the proposal equation,
\begin{equation*}
    \norm{x^\star -x} = \norm{-\delta \nabla U(x_{\frac{1}{2}})+\sqrt{2\delta}\xi} \leq \delta\norm{\nabla U(x_{\frac{1}{2}})}+\sqrt{2\delta}\norm{\xi}.
\end{equation*}
Applying~\eqref{eq:first_bound},
\begin{equation*}
    \norm{x^\star -x} \leq \delta L\norm{x_{\frac{1}{2}}}+\sqrt{2\delta}\norm{\xi}.
\end{equation*}
Therefore, since $|s|\leq 1/2$,
\begin{align*}
    \norm{y_s}
    &\geq \left(1-\frac{L\delta}{2}\right)\norm{x_{\frac{1}{2}}}
    -\sqrt{\frac{\delta}{2}}\norm{\xi}.
\end{align*}
If $L\delta<2$, a two-case argument, according to whether $\norm{x_{\frac{1}{2}}}$ is larger or smaller than a constant multiple of $\norm{\xi}$, shows that
\begin{equation*}
    \frac{\norm{x_{\frac{1}{2}}}^3}{(1+\norm{y_s})^3}
    \leq C_\delta\left(1+\norm{\xi}^3\right),
\end{equation*}
uniformly in $s\in[-1/2,1/2]$ and $x$. The same argument applies to $A_1$, $A_2$, and $A_3$, and therefore
\begin{equation*}
    |\Delta(x,x^\star)|\leq C_\delta\left(1+\norm{\xi}^3\right).
\end{equation*}
Thus, cubic decay of $\nabla^3U$ bounds the energy error uniformly by an integrable random variable. Moreover,
\begin{equation*}
    \inf_{x\in\mathbb{R}^d}\mathbb{E}[\alpha_\delta(x,x^\star)]
    \geq \mathbb{E}\left[\exp\left(-C_\delta(1+\norm{\xi}^3)\right)\right]
    >0.
\end{equation*}
Thus, Assumption~\ref{ass:third_derivative_decay} provides a uniform positive lower bound on the mean acceptance probability.
\end{proof}

\begin{remark}
    To illustrate the role of the decay-rate condition in Assumption~\ref{ass:third_derivative_decay}, we consider a one-dimensional version of the Bayesian inverse-problem toy model in \cite{DM15}. The observation model is
\begin{equation*}
    b=x+\varepsilon,
\end{equation*}
where $b\in\mathbb{R}$ is the observed data, $x\in\mathbb{R}$ is the unknown signal, and $\varepsilon\sim\mathcal{N}(0,1)$ is the observation noise. We use a smoothed exponential-power prior. Up to an additive constant, the negative log-posterior is
\begin{equation*}
    U(x) = \frac{1}{2}(x-b)^2 + (x^2+\rho)^\beta,
\end{equation*}
where $\rho>0$ and $\beta\in(1/2,1)$.

Since $(x^2+\rho)^\beta$ is convex for $\beta\geq 1/2$, $U$ is $1$-strongly convex. Moreover, $U'$ is $L$-Lipschitz, with
\begin{equation*}
    L=1+\frac{2\beta}{\rho^{1-\beta}}.
\end{equation*}
The third derivative is
\begin{align*}
    U'''(x) &=-12\beta (1-\beta)\frac{x}{(x^2+\rho)^{2-\beta}}
    +8\beta (1-\beta)(2-\beta)\frac{x^3}{(x^2+\rho)^{3-\beta}}.
\end{align*}
As $|x|\to\infty$,
\begin{equation*}
    U'''(x)=-4\beta(1-\beta)(2\beta-1)\,\operatorname{sign}(x)|x|^{2\beta-3}
    +o\bigl(|x|^{2\beta-3}\bigr).
\end{equation*}
Thus, $|U'''(x)|(1+|x|)^3$ grows like $|x|^{2\beta}$ and is unbounded. In particular, $U'''$ decays, but not at the rate required by Assumption~\ref{ass:third_derivative_decay}.

Figure~\ref{fig:inverse_problem_energy_beta} shows Monte Carlo estimates of the mean energy error, based on $5\times10^4$ independent proposals for each value of $x$. We use $\delta=0.1$, $b=0$, and $\rho=10^{-4}$. Figure~\ref{fig:inverse_problem_acceptance_beta} shows the corresponding estimates of the mean acceptance probability, which decrease toward zero over the plotted range as $|x|$ increases.

\begin{figure}[ht]
    \centering
    \begin{subfigure}{0.48\textwidth}
        \centering
        \includegraphics[width=\linewidth]{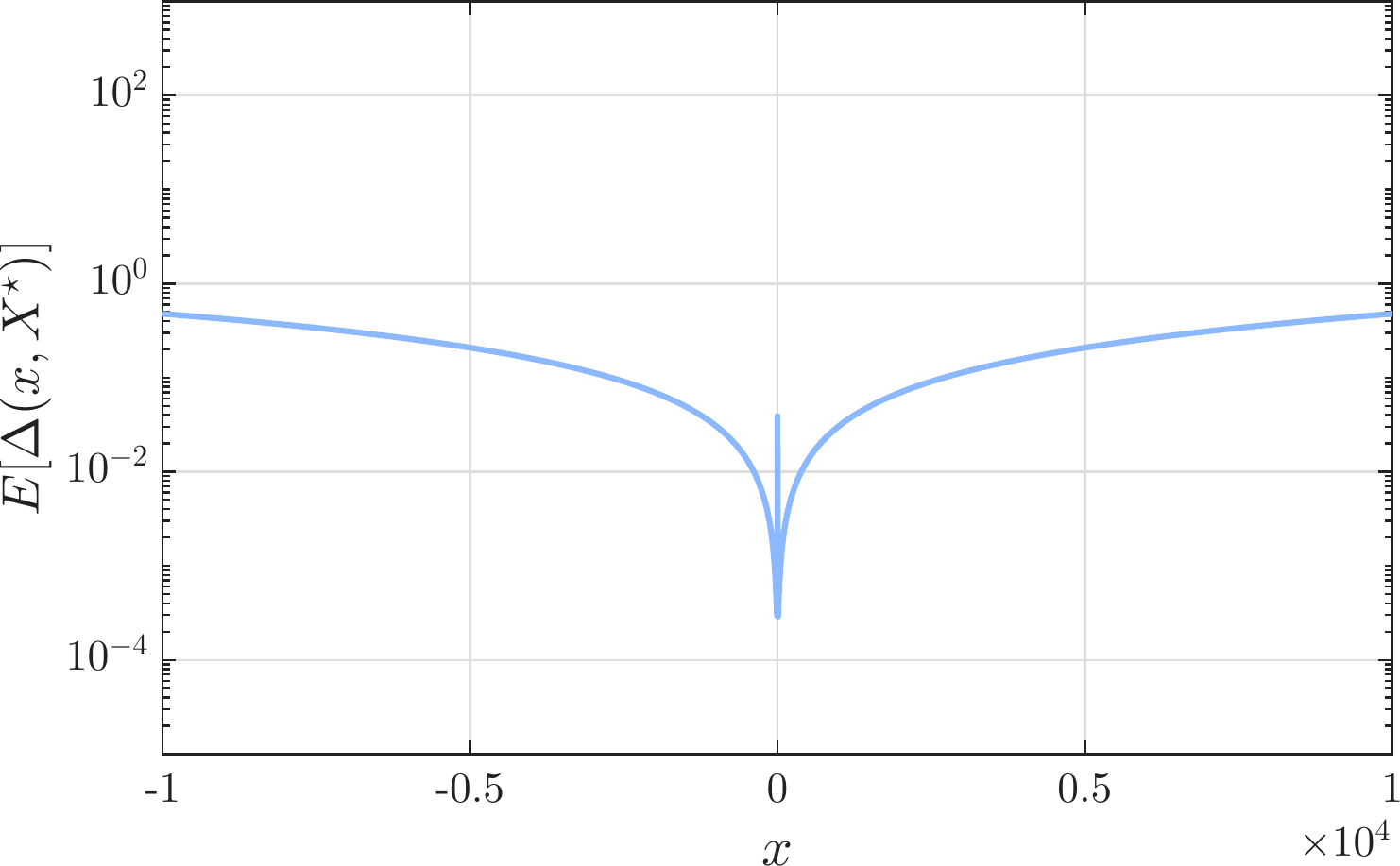}
        \subcaption{$\mathbb{E}[\Delta(x,x^\star)]$ for $\beta=0.6$.}
    \end{subfigure}\hfill
    \begin{subfigure}{0.48\textwidth}
        \centering
        \includegraphics[width=\linewidth]{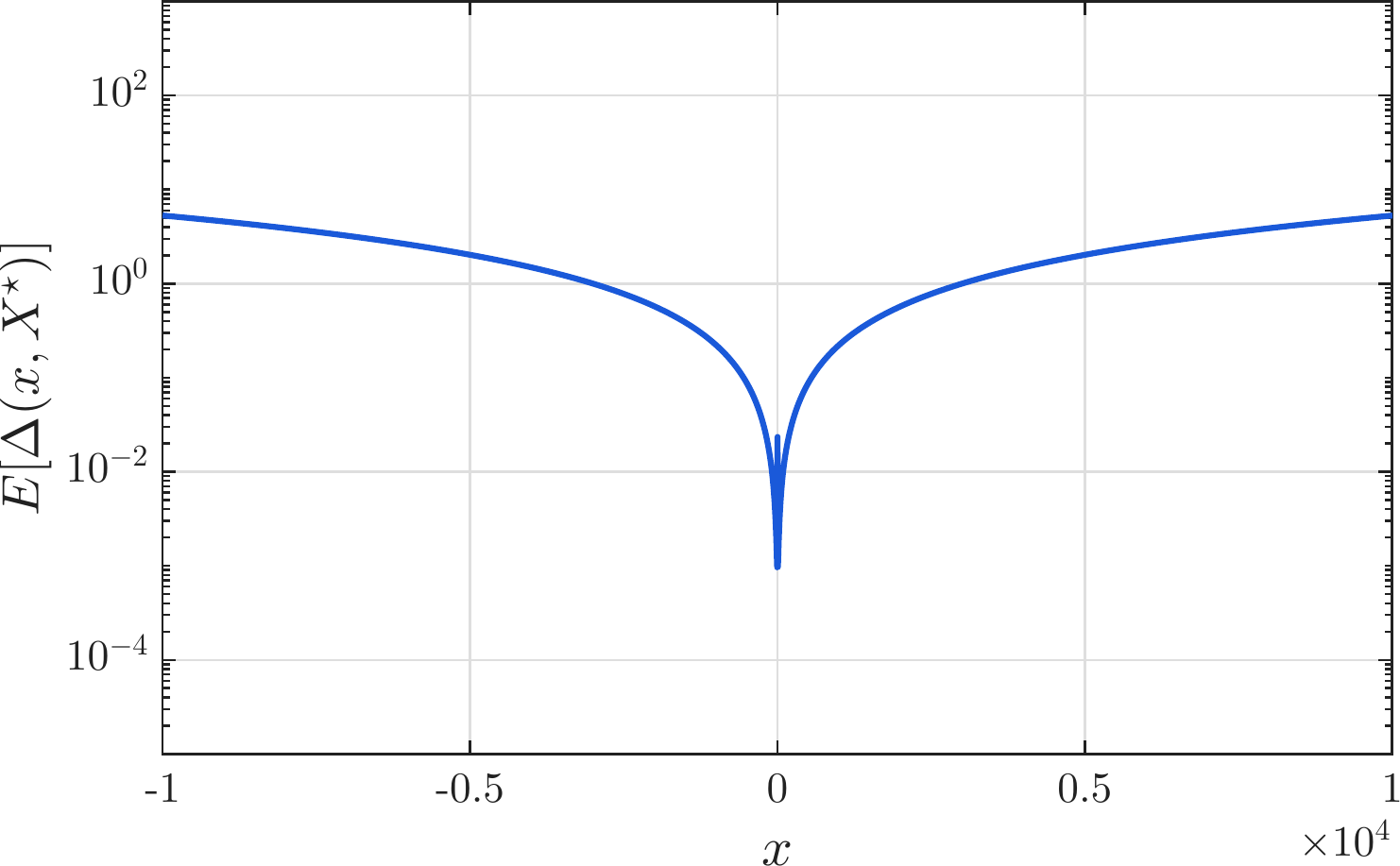}
        \subcaption{$\mathbb{E}[\Delta(x,x^\star)]$ for $\beta=0.7$.}
    \end{subfigure}

    \vspace{0.8em}

    \begin{subfigure}{0.48\textwidth}
        \centering
        \includegraphics[width=\linewidth]{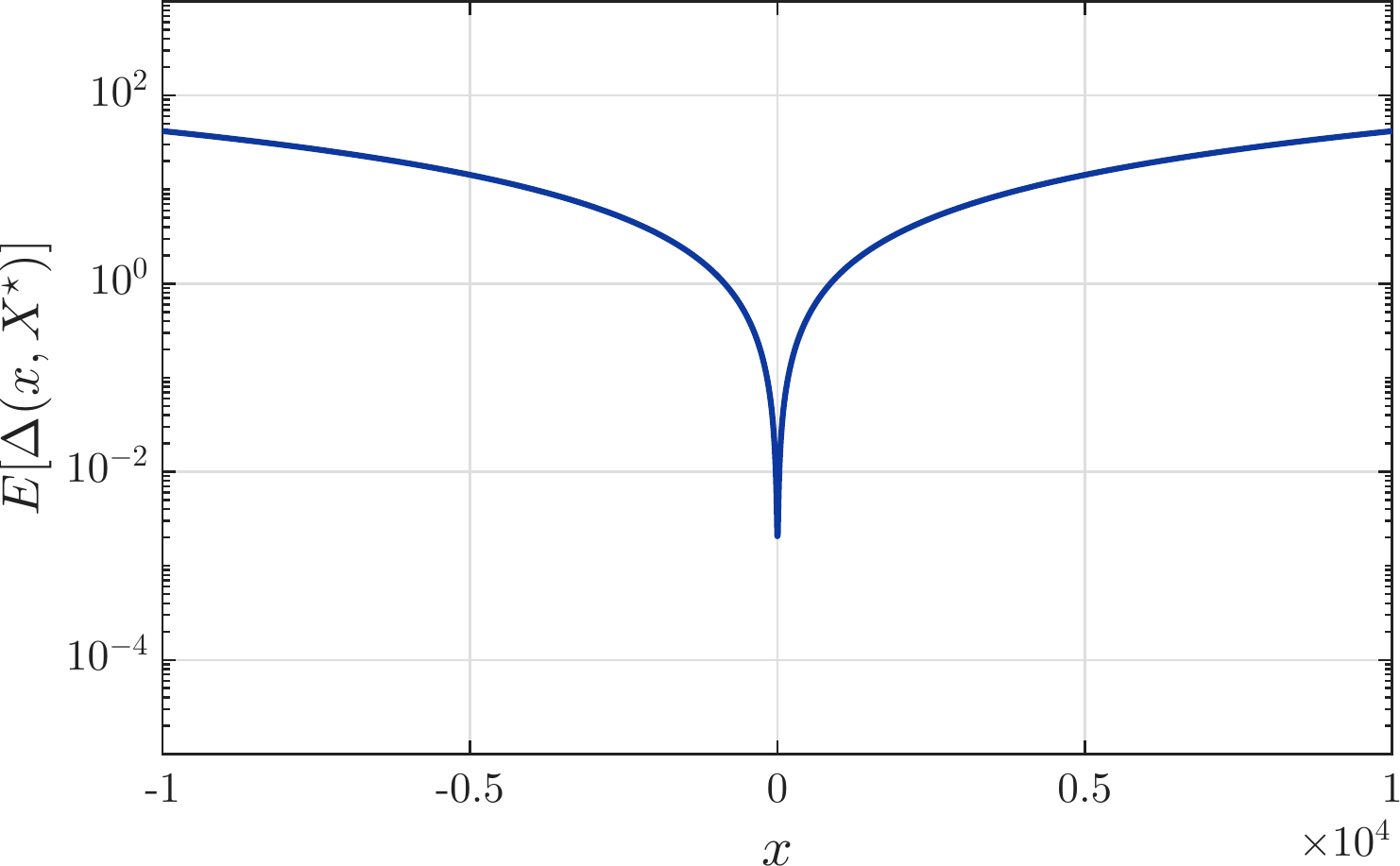}
        \subcaption{$\mathbb{E}[\Delta(x,x^\star)]$ for $\beta=0.8$.}
    \end{subfigure}\hfill
    \begin{subfigure}{0.48\textwidth}
        \centering
        \includegraphics[width=\linewidth]{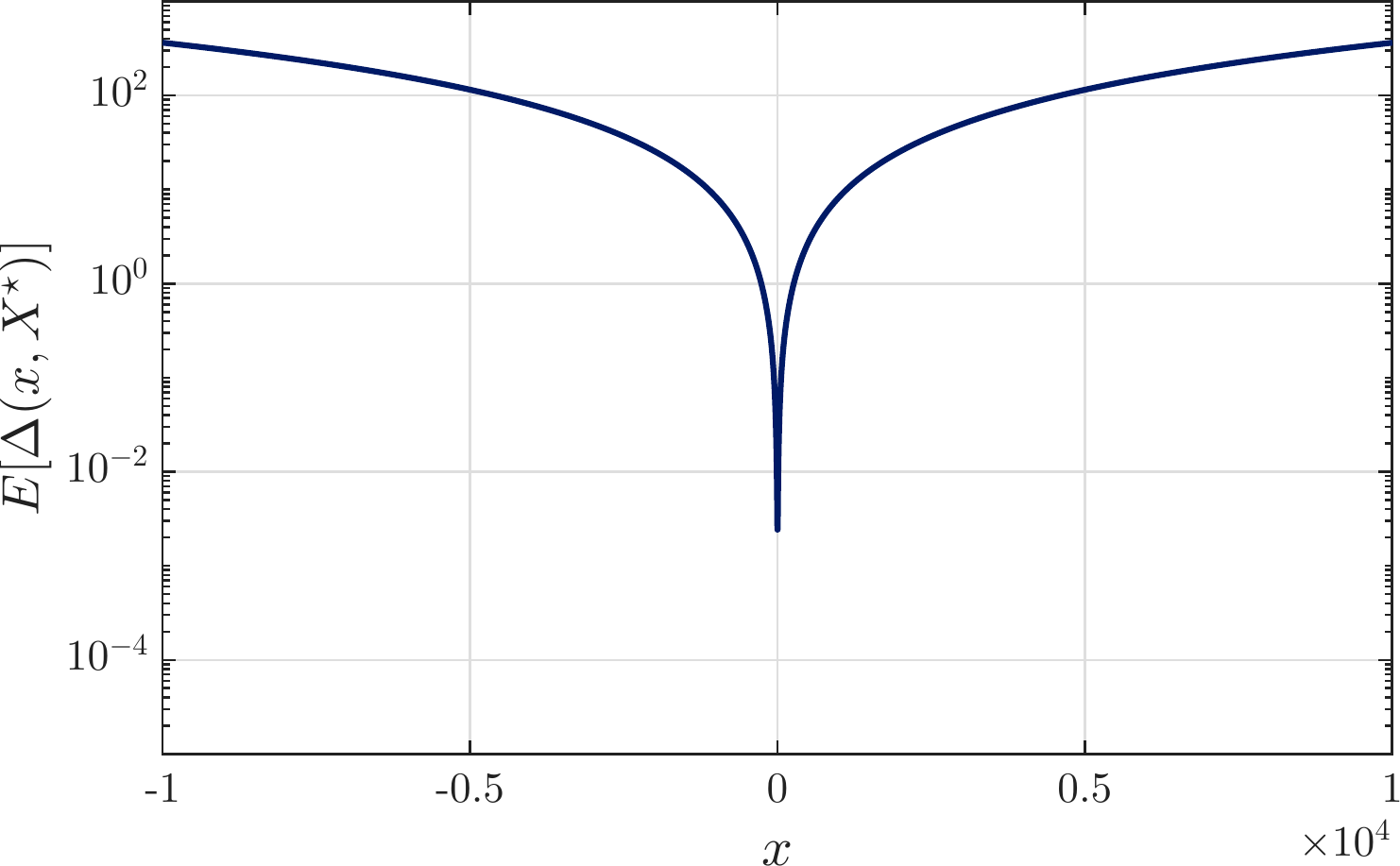}
        \subcaption{$\mathbb{E}[\Delta(x,x^\star)]$ for $\beta=0.9$.}
    \end{subfigure}

    \caption{Expected energy error $\mathbb{E}[\Delta(x,x^\star)]$ as a function of $x$ for different values of $\beta$.}
    \label{fig:inverse_problem_energy_beta}
\end{figure}

\begin{figure}[ht]
    \centering
    \includegraphics[width=0.6\linewidth]{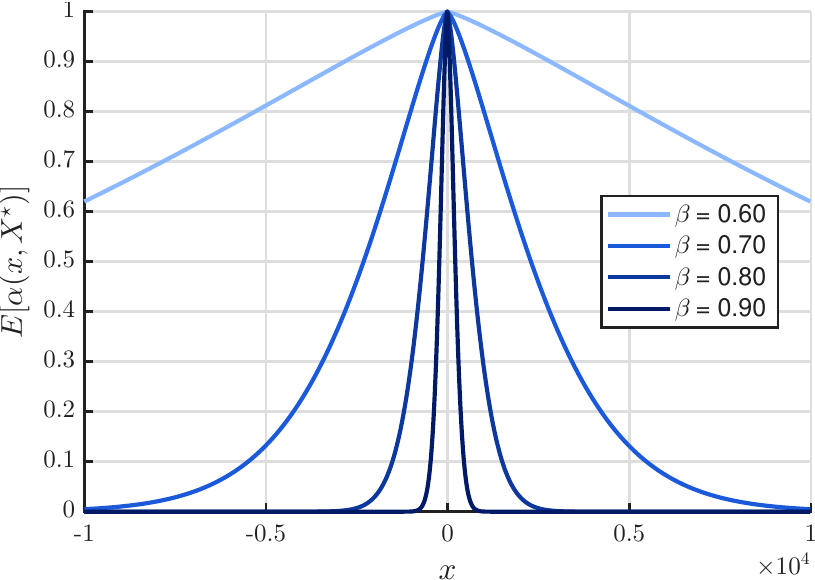}

    \caption{Expected acceptance probability $\mathbb{E}[\alpha_\delta(x,x^\star)]$ as a function of $x$ for different values of $\beta$.}
    \label{fig:inverse_problem_acceptance_beta}
\end{figure}

\end{remark}

\section{Numerical investigations}\label{sec:num}

\subsection{Gaussian mixture}\label{sec:gaussian_mixture}

We first consider the $d$-dimensional Gaussian-mixture target
\begin{equation*}
    \pi(x)
    \propto
    \exp\left(-\frac12\norm{x-a}^2\right)
    +\exp\left(-\frac12\norm{x+a}^2\right),
    \qquad a=c e_1,
\end{equation*}
as studied in \cite{Dalalyan2017,dwivedi18a}. Up to an additive constant, its potential is
\begin{equation*}
    U(x)=\frac12\norm{x}^2-\log\!\left(\cosh(a^{\top}x)\right).
\end{equation*}
Unlike the example considered in Section~\ref{sec:optimal_scaling}, this target is not an i.i.d.\ product distribution. However,  $\pi$ is a product of a one dimensional Gaussian mixture with $d-1$ standard Gaussian distributions.  Consequently, only the first coordinate contributes to the MAIMLA acceptance probability.

\begin{figure}[t!]
    \centering
    \includegraphics[width=0.8\textwidth]{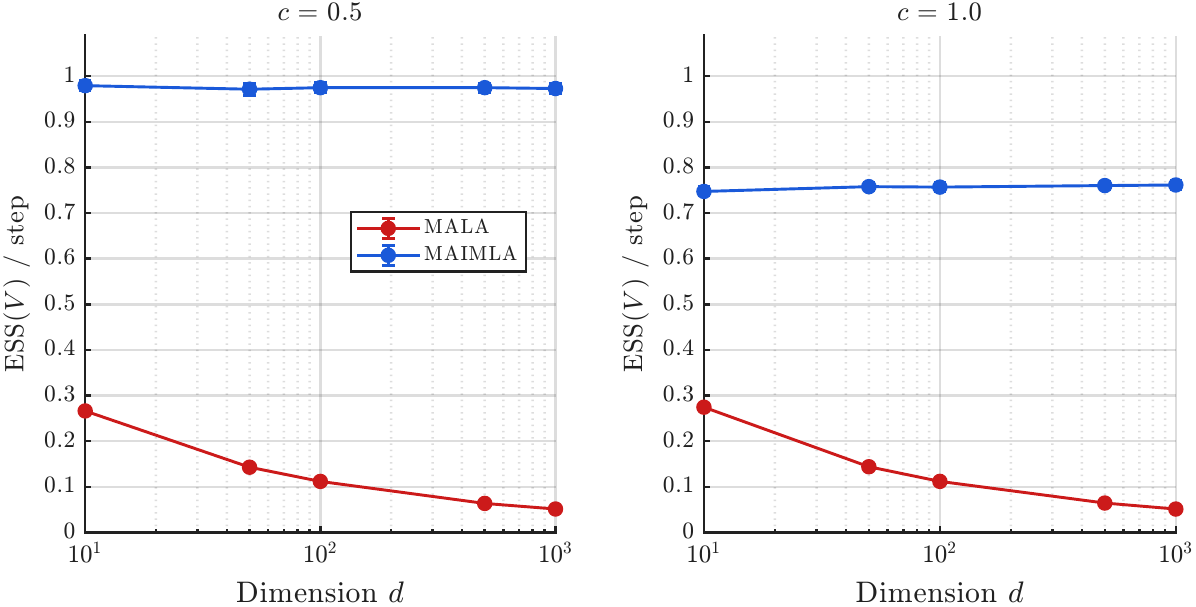}
    \caption{Effective sample size per iteration for the $d$-dimensional Gaussian-mixture target, averaged over $50$ runs with $N=10^5$ samples per run. Error bars show one empirical standard deviation.}
    \label{fig:ess_gauss_mix}
\end{figure}

\begin{figure}[ht!]
    \centering
    \includegraphics[width=0.8\textwidth]{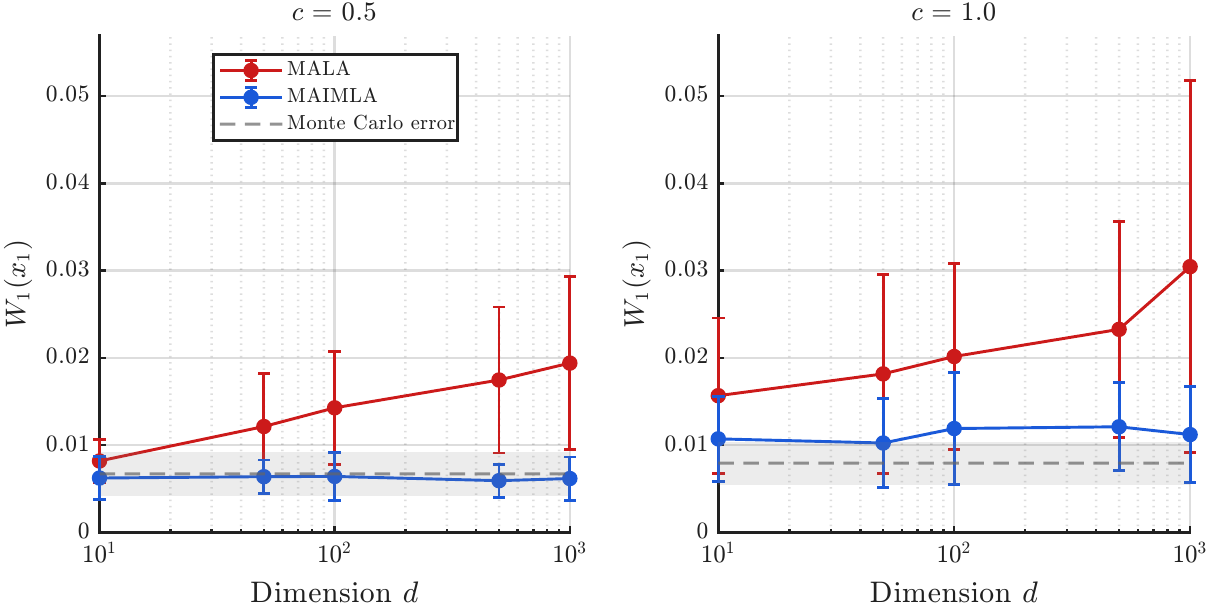}
    \caption{Wasserstein-$1$ error for the $d$-dimensional Gaussian-mixture target, averaged over $50$ runs with $N=10^5$ samples per run. The shaded grey region represents the reference Monte Carlo error.}
    \label{fig:w1_gauss_mix}
\end{figure}

We compare MALA and MAIMLA for $c\in\{0.5,1\}$ and dimensions
$d\in\{10,50,100,500,1000\}$. For $c<1$, the potential is $L$-gradient Lipschitz and $m$-strongly convex with   $L=1$ and  $m=1-c^2$. In this setting MAIMLA accepts with high probability for all step sizes so it is not possible to achieve an average acceptance probability of 0.574. Instead  we use the step size $\delta=\frac{2}{\sqrt{Lm}}$ which achieves the fastest contraction rate for the unadjusted algorithm \cite{KDAPSZ24}. When $c=1$, the potential is convex but not strongly convex, and we instead set $\delta=1$. The MALA step size is tuned to give an average acceptance probability of approximately $0.574$. For each dimension and each value of $c$, we generate $N=10^5$ samples and repeat the experiment $50$ times. We report the effective sample size per iteration for $U(x)$ in Figure \ref{fig:ess_gauss_mix} and the Wasserstein-$1$ error of the first component $x_{1}$ in Figure \ref{fig:w1_gauss_mix}, together with their empirical standard deviations across the repeated runs. In both metrics, the performance of  MALA deteriorates as the dimension $d$ increases while the performance of MAIMLA remains constant across all dimensions. This is to be expected since MAIMLA is capable of operating with a dimension independent time-step unlike MALA.

\subsection{Image denoising}\label{sec:image_denoising}

We next consider an image-denoising inverse problem. Let $x\in\mathbb{R}^d$ denote the unknown vectorised image and suppose that the observation is
\begin{equation*}
    y=x+\xi,
    \qquad
    \xi\sim\mathcal{N}(0,\sigma^2I_d).
\end{equation*}
The corresponding likelihood is
\begin{equation*}
    p(y\mid x)
    \propto
    \exp\left(-\frac{\norm{x-y}^2}{2\sigma^2}\right).
\end{equation*}

To regularise this ill-posed problem we  impose, similarly to   \cite{ASS98, CPM18},  the  following smooth sparsity-promoting prior
\begin{equation*}
    p(x)
    \propto
    \exp\left(
        -\lambda\sum_{j=1}^{d_{\mathrm{det}}}
        H_{\varepsilon}\bigl((W_{\mathrm{det}}x)_j\bigr)
    \right),
\end{equation*}
where $W_{\mathrm{det}}$ contains the detail rows of an orthonormal wavelet transform and
\begin{equation*}
    H_{\varepsilon}(t)
    =
    \begin{cases}
        \dfrac{t^2}{2\varepsilon}, & |t|\leq\varepsilon,\\[6pt]
        |t|-\dfrac{\varepsilon}{2}, & |t|>\varepsilon
    \end{cases}
\end{equation*}
is the Huber function. The posterior potential is therefore
\begin{equation*}
    U(x)
    =\frac{\norm{x-y}^2}{2\sigma^2}
    +\lambda\sum_{j=1}^{d_{\mathrm{det}}}
        H_{\varepsilon}\bigl((W_{\mathrm{det}}x)_j\bigr).
\end{equation*}
Observe that $U$ is convex and has $L$-Lipschitz gradient with  $L=1/\sigma^2 + \lambda /\epsilon$.
 Increasing $\varepsilon$ enlarges the quadratic region of  the Huber penalty and makes the posterior increasingly Gaussian-like, providing a useful regime in which to compare MALA and MAIMLA. Moreover, because the likelihood is quadratic, the wavelet transform is orthonormal, and the Huber penalty is separable in the wavelet coefficients, the implicit midpoint equation decouples coefficientwise after transformation to the wavelet domain. Its solution can therefore be evaluated in closed form, without requiring a generic nonlinear solver at each MAIMLA iteration.  We note here that  when $\varepsilon=0$, $p(x)$ is a  Laplace-type wavelet prior, hence the   potential is nonsmooth and therefore lies outside the setting considered here.

For our numerical experiment we use  the cameraman image ($256 \times 256$ pixels) with $\sigma=0.1$,   $\varepsilon\in\{0.5\sigma,1\sigma,2\sigma\} = \{0.05, 0.1, 0.2\}$, and compare MALA and MAIMLA.  In both cases the chains are initialised at $X_0 = y$ with $10^4$ burn-in iterations and  $5 \times 10^5$ sampling iterations. The step size is chosen adaptively using the Robbins--Monro rule \cite{YAJR05} targeting 0.574 acceptance probability with initial step size $\delta=10^{-2} \times  \delta_{\text{ref}}$, where $\delta_{\text{ref}} = 0.9 \times 2/L$. We record summary statistics such as the posterior mean for each sampler,
and traces of 300 randomly chosen pixels from which we compute autocorrelation and effective sample size (ESS). Alongside mean square error (MSE) \footnote{$\text{MSE} = (1/n) \sum_{i=1}^n (\hat{x}_i - x_i)^2$ with $\hat{x}$ the posterior mean and $x$ the ground truth, $n$ the number of pixels}, we report
structural similarity index measure (SSIM)\footnote{Higher is better, $\text{SSIM}=1$ means the images are identical.}\cite{SSIM04} as image evaluation metrics.

\begin{figure}[p]
    \centering
    \begin{minipage}[t]{0.31\textwidth}
        \centering
        \includegraphics[width=\linewidth]{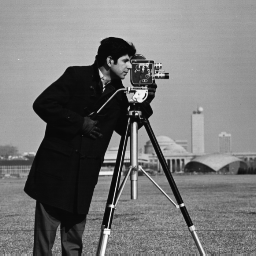}\\[-0.4em]
        {\scriptsize Ground truth}
    \end{minipage}\hfill
    \begin{minipage}[t]{0.31\textwidth}
        \centering
        \includegraphics[width=\linewidth]{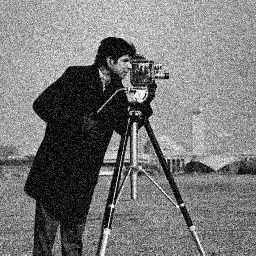}\\[-0.4em]
        {\scriptsize Noisy observation}
    \end{minipage}\hfill
    \begin{minipage}[t]{0.31\textwidth}
        \centering
        \includegraphics[width=\linewidth]{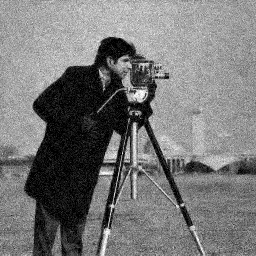}\\[-0.4em]
        {\scriptsize MAIMLA posterior mean}
    \end{minipage}
    \caption{Ground truth, noisy observation, and MAIMLA posterior mean reconstruction for the denoising experiment with Huber width $\varepsilon=0.1$. The corresponding MALA reconstruction is visually indistinguishable from the MAIMLA reconstruction. Reconstructions at the other tested values of $\varepsilon$ differ only slightly and are omitted for clarity.}
    \label{fig:denoising_mmse}
\end{figure}

\begin{figure}[p]
    \centering
    \begin{minipage}[t]{0.31\textwidth}
        \centering
        \textbf{$\varepsilon=0.05$}\par
        \includegraphics[width=\linewidth]{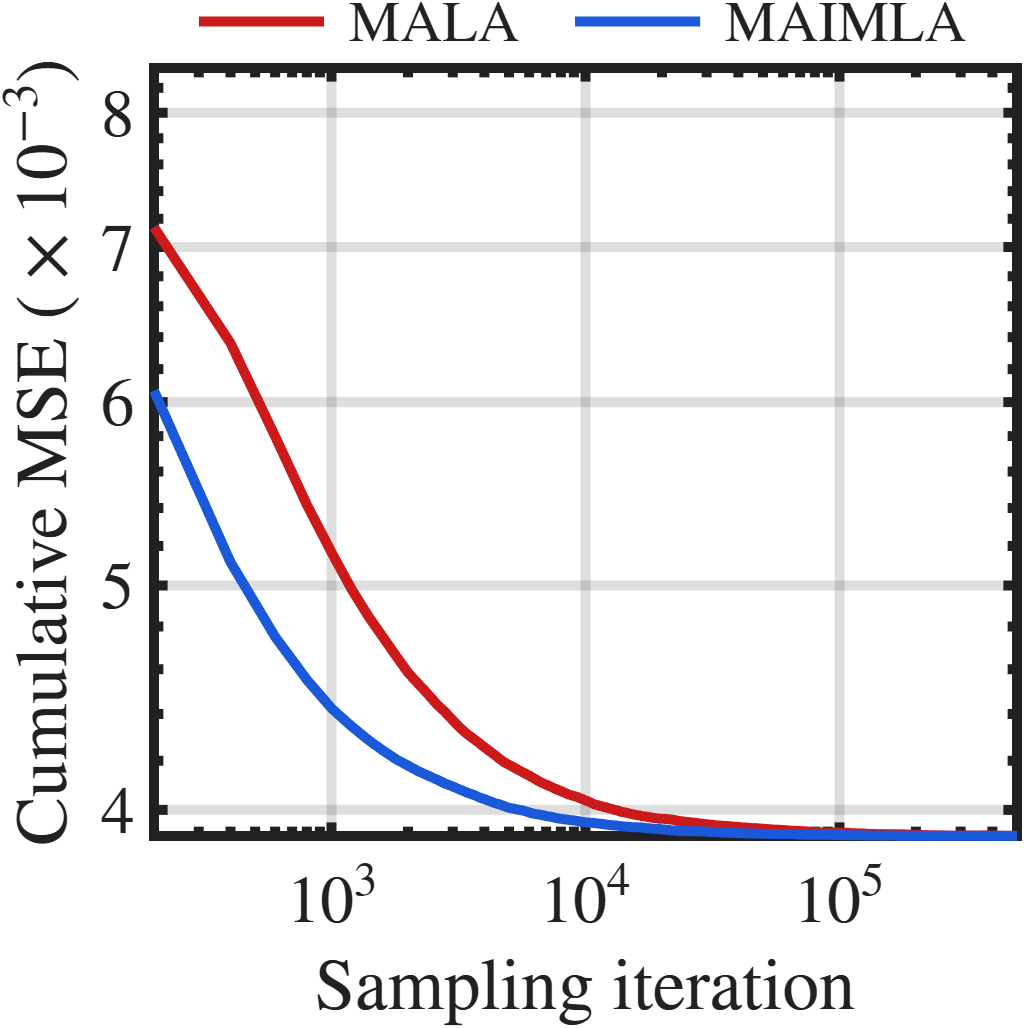}
    \end{minipage}\hfill
    \begin{minipage}[t]{0.31\textwidth}
        \centering
        \textbf{$\varepsilon=0.10$}\par
        \includegraphics[width=\linewidth]{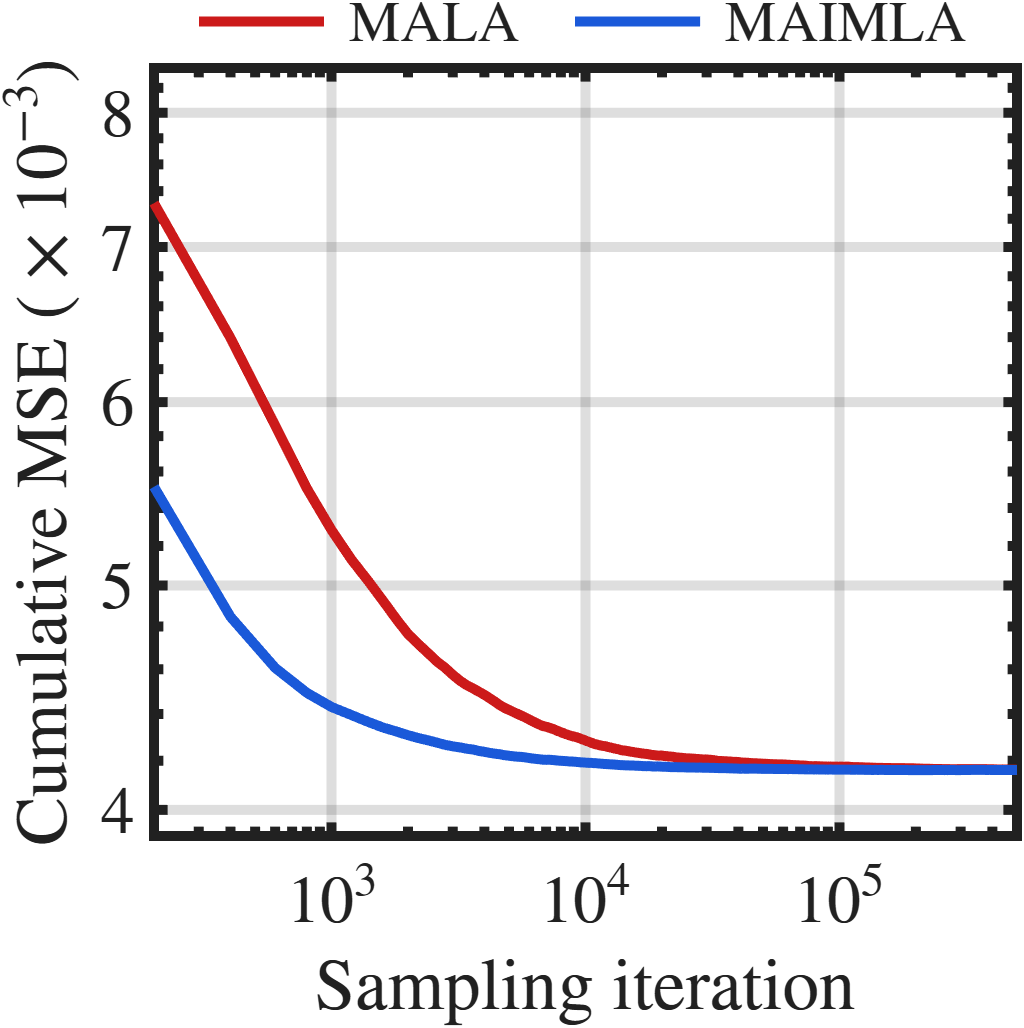}
    \end{minipage}\hfill
    \begin{minipage}[t]{0.31\textwidth}
        \centering
        \textbf{$\varepsilon=0.20$}\par
        \includegraphics[width=\linewidth]{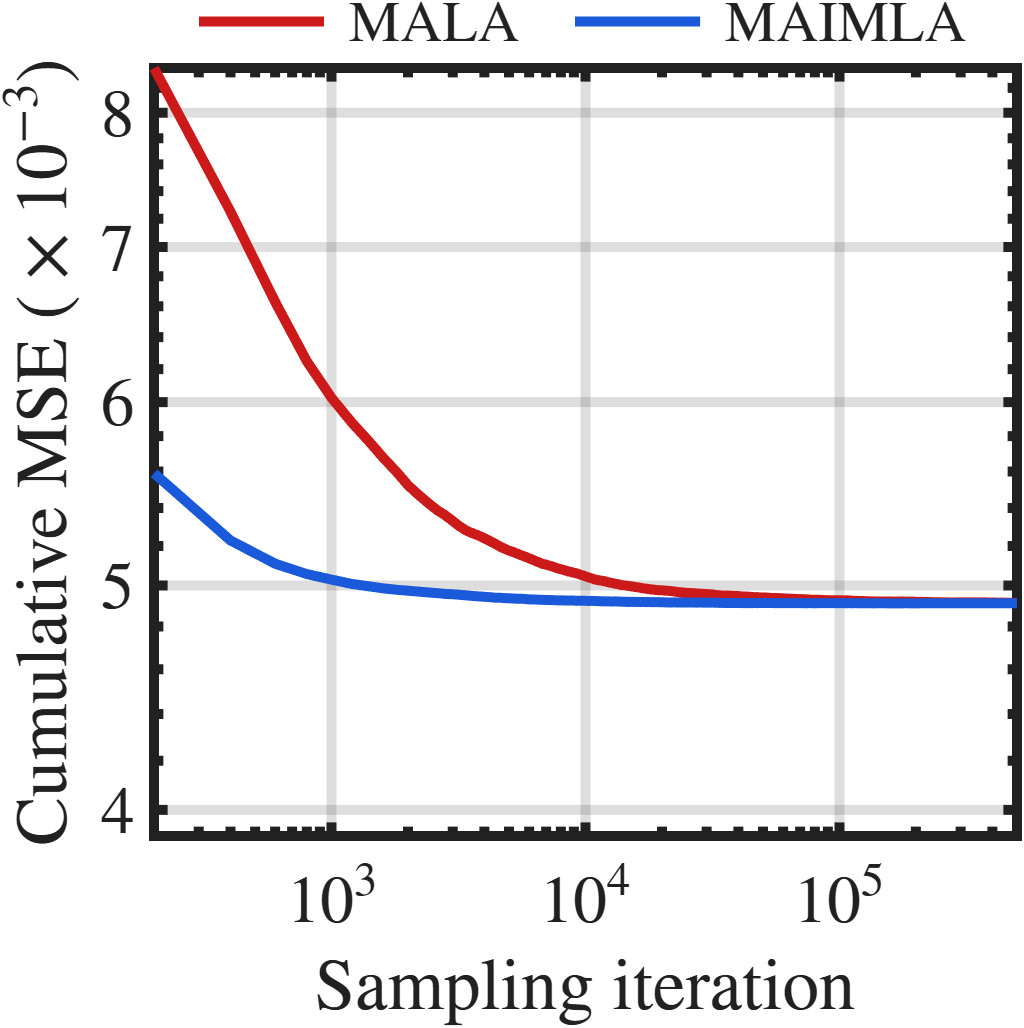}
    \end{minipage}
    \caption{Cumulative mean squared error (MSE) for the denoising experiment at the three Huber widths $\epsilon$. Each panel compares MALA and MAIMLA.}
    \label{fig:denoising_mse_cumulative}
\end{figure}

\begin{figure}[p]
    \centering
    \begin{minipage}[t]{0.31\textwidth}
        \centering
        $\varepsilon=0.05$\par
        \includegraphics[width=\linewidth]{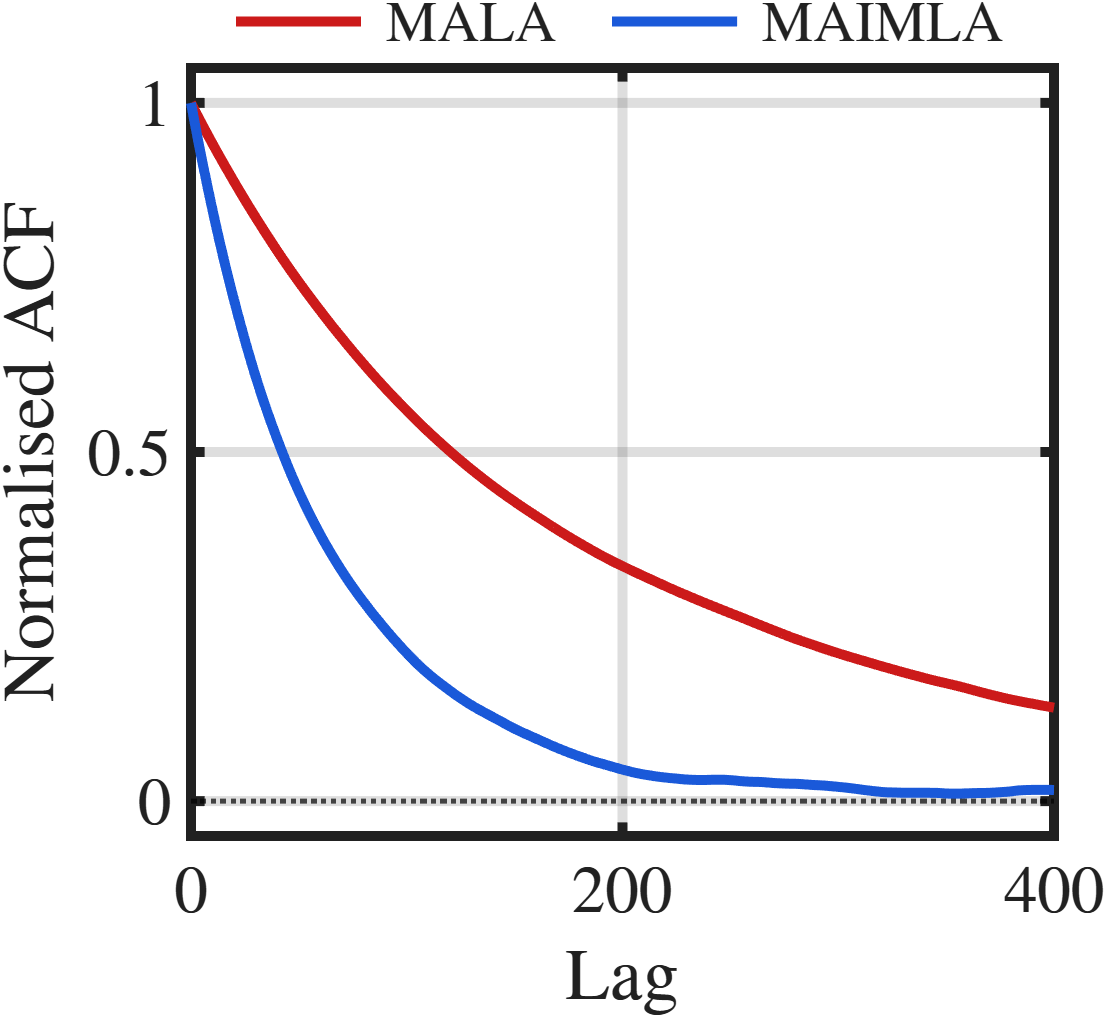}\\[-0.4em]
        {\scriptsize Slowest pixel}\\[0.4em]
        \includegraphics[width=\linewidth]{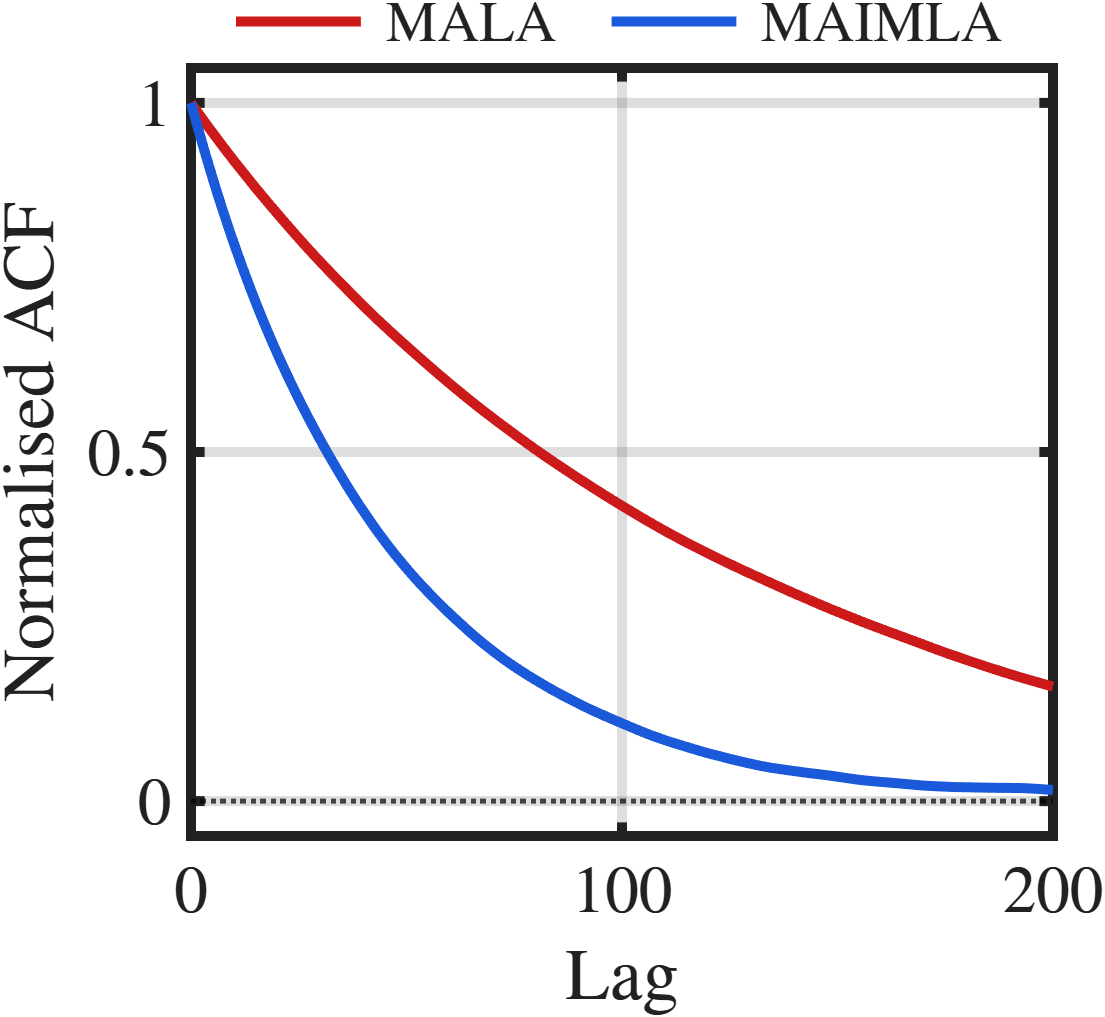}\\[-0.4em]
        {\scriptsize Fastest pixel}
    \end{minipage}\hfill
    \begin{minipage}[t]{0.31\textwidth}
        \centering
        \textbf{$\varepsilon=0.10$}\par
        \includegraphics[width=\linewidth]{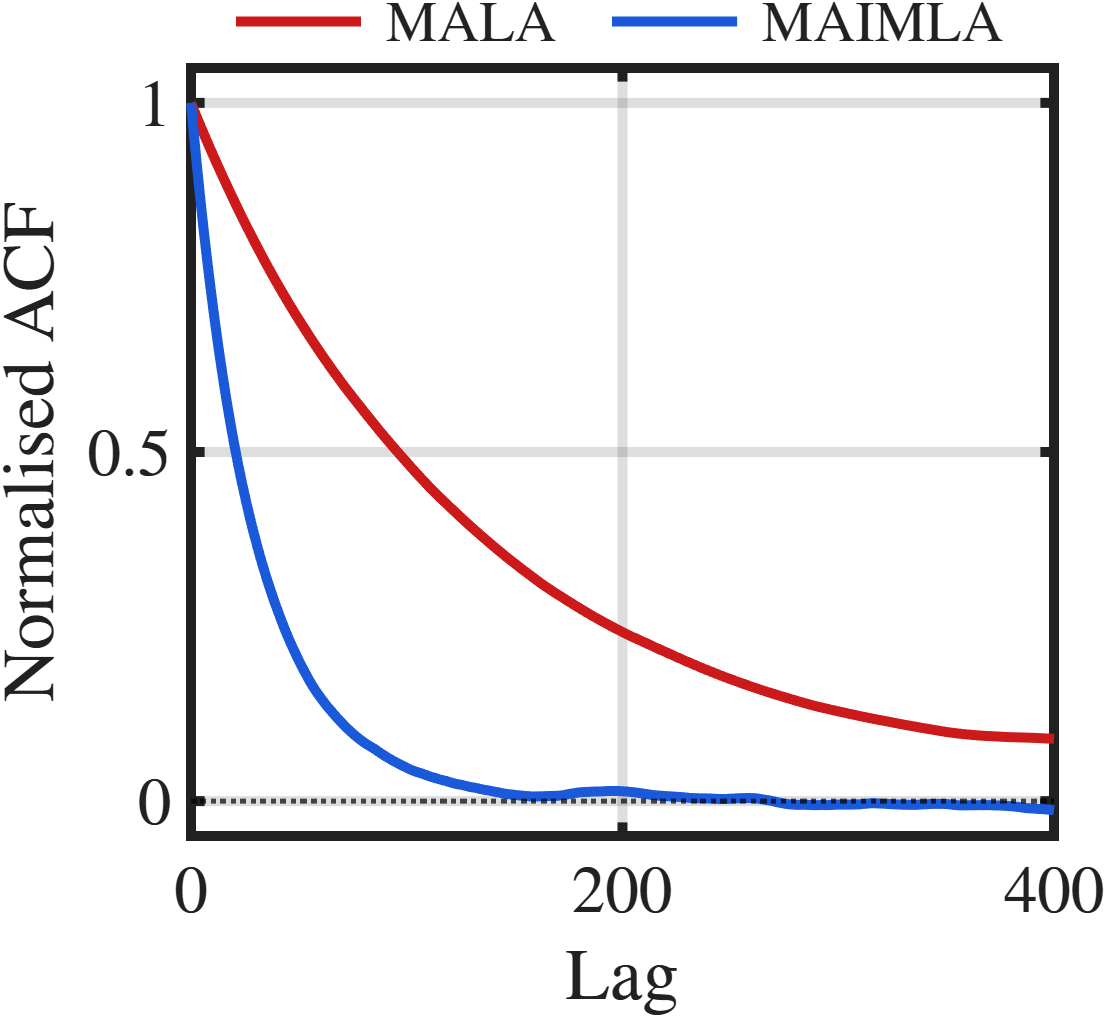}\\[-0.4em]
        {\scriptsize Slowest pixel}\\[0.4em]
        \includegraphics[width=\linewidth]{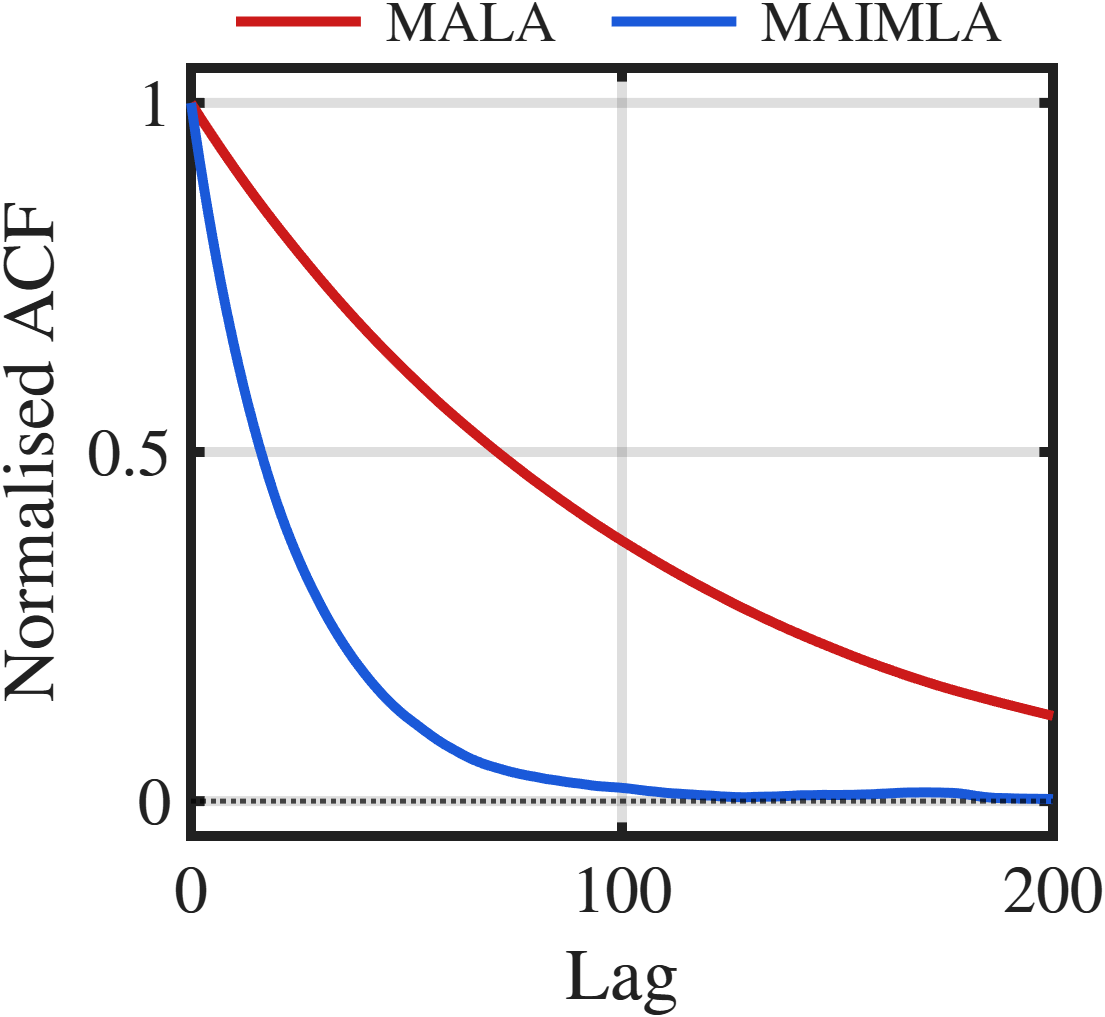}\\[-0.4em]
        {\scriptsize Fastest pixel}
    \end{minipage}\hfill
    \begin{minipage}[t]{0.31\textwidth}
        \centering
        \textbf{$\varepsilon=0.20$}\par
        \includegraphics[width=\linewidth]{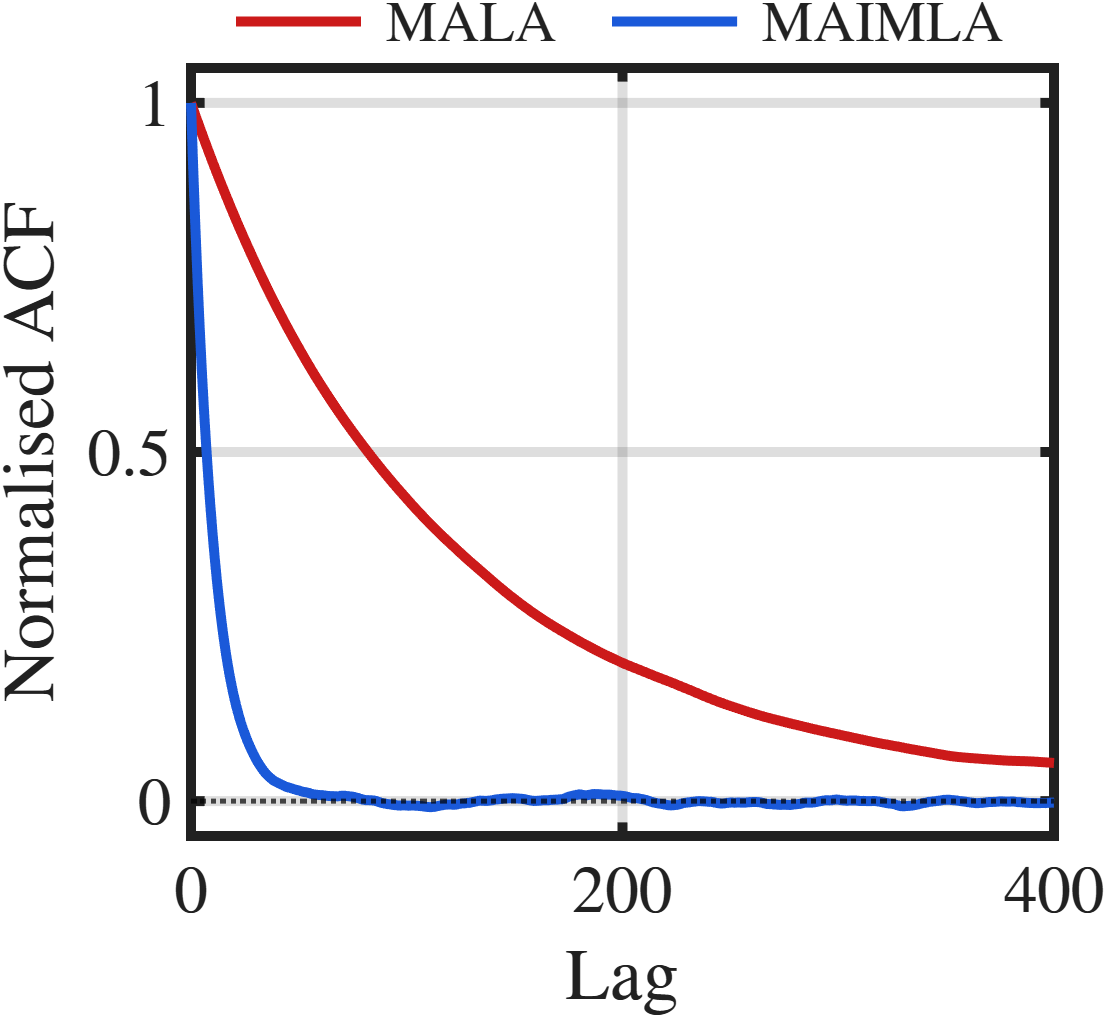}\\[-0.4em]
        {\scriptsize Slowest pixel}\\[0.4em]
        \includegraphics[width=\linewidth]{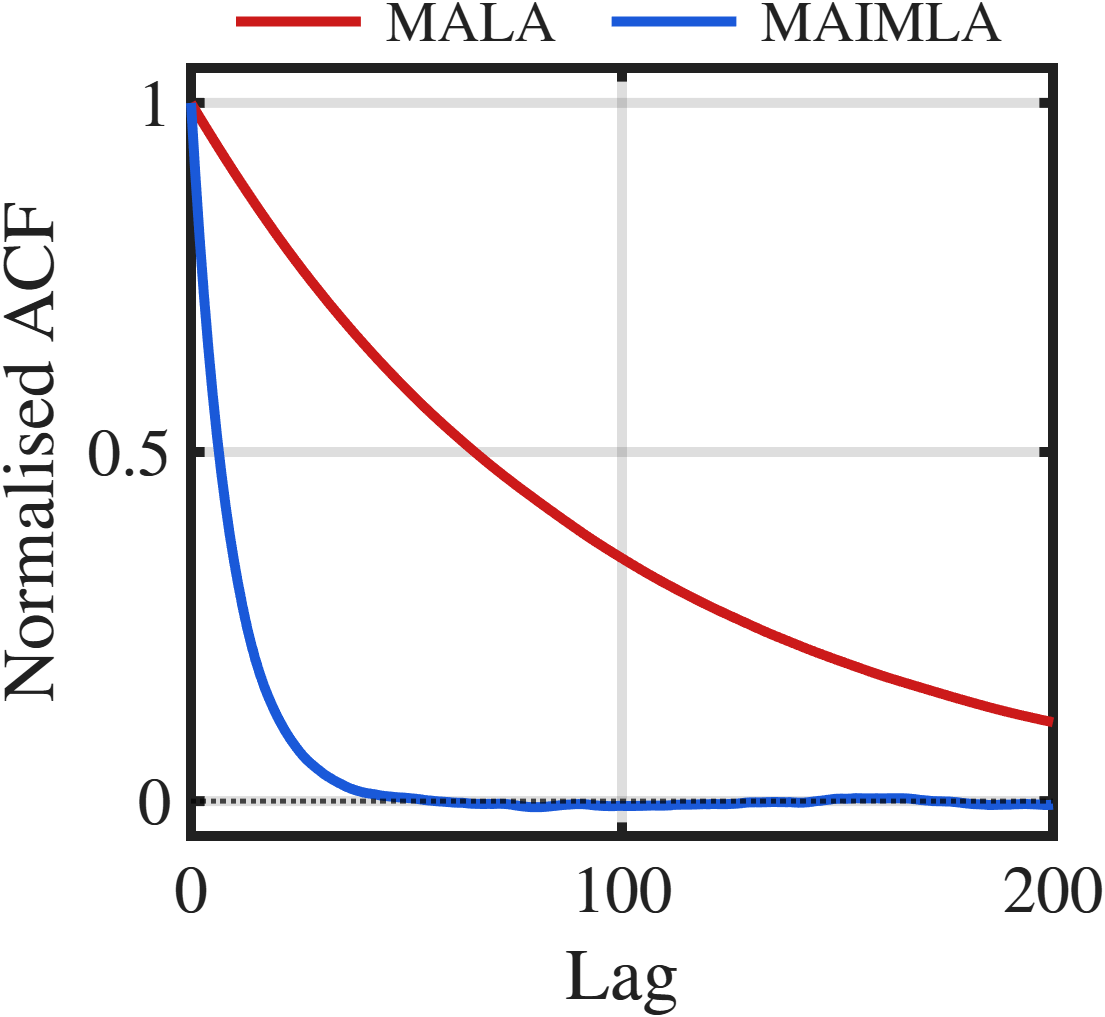}\\[-0.4em]
        {\scriptsize Fastest pixel}
    \end{minipage}
    \caption{Autocorrelation functions for the slowest- and fastest-mixing pixels in the denoising experiment. Columns correspond to the Huber width $\varepsilon$; the top and bottom rows show the slowest and fastest pixels, respectively.}
    \label{fig:denoising_acf}
\end{figure}

Table~\ref{tab:denoising_sampler_results} summarises the mean acceptance probabilities, effective step sizes, posterior mean reconstruction quality, and ESS for the two samplers. The ESS column reports the minimum and maximum across the 300 tracked pixels.
\begin{table}[htbp]
{
    \centering
    \scriptsize
    \setlength{\tabcolsep}{3pt}
    \begin{tabular}{@{}lllrrrr@{}}
        \toprule
        $\varepsilon$ & Sampler & Acc. & $\delta$ ($10^{-4}$) & MSE ($10^{-3}$) $\downarrow$ & SSIM $\uparrow$ & ESS range \\
        \midrule
        0.05 & MALA   & 0.56 & 0.769 & 3.90 & 0.4662 &  1532.5--2319.4 \\
        0.05 & MAIMLA & 0.59 & 1.93 & 3.90 & 0.4663 &  3838.0--5372.8 \\
        0.10 & MALA   & 0.56 & 0.975 & 4.16 & 0.4522 &  1869.4--2683.0 \\
        0.10 & MAIMLA & 0.58 & 4.30 & 4.16 & 0.4522 &  7822.8--10352.5 \\
        0.20 & MALA   & 0.58 & 1.12 & 4.91 & 0.4213 &  2118.8--2764.8 \\
        0.20 & MAIMLA & 0.57 & 13.0 & 4.91 & 0.4214 & 23292.1--26715.5 \\
        \bottomrule
    \end{tabular}
    \caption{Posterior mean reconstruction results for MALA and MAIMLA, including mean acceptance probabilities, step sizes $\delta$, MSE, SSIM, and ESS range over 300 tracked pixels. The $\delta$ and MSE columns are reported in units of $10^{-4}$ and $10^{-3}$, respectively.}
    \label{tab:denoising_sampler_results}
}
\end{table}
As we can see in Table \ref{tab:denoising_sampler_results}  both algorithms
have the same reconstruction quality, as reported in terms of MSE and SSIM of the posterior mean for each $\epsilon$ (see also Figure \ref{fig:denoising_mmse}), which indicates that they   are both close to stationarity. However, MAIMLA converges faster to stationarity as illustrated in Figure \ref{fig:denoising_mse_cumulative} which shows the convergence of the MSE as a function of the iteration number. A similar pattern is observed in Figure \ref{fig:denoising_acf} which plots the normalised autocorrelation for the slowest and fastest pixel. The reported slowest-pixel ESS of MAIMLA is about 2.5, 4.2, and 11 times MALA’s at $\epsilon = \{0.05,0.1,0.2\}$, respectively. This is consistent with the fact that MAIMLA can achieve the same acceptance probability as MALA, but with a larger timestep.

\section{Conclusions} \label{sec:con}

We have studied the implicit midpoint Langevin algorithm IMLA and its Metropolis-adjusted variant, MAIMLA. The analysis establishes geometric ergodicity of IMLA under convexity, dissipativity, and polynomial growth conditions on the gradient, including the case where the gradient grows superlinearly. For MAIMLA, geometric ergodicity is established both for convex potentials with a radially symmetric leading-order term in the superlinear-growth regime and, under a suitable third-derivative decay condition, when the gradient is globally Lipschitz. A central consequence is that MAIMLA can remain geometrically ergodic in superlinear settings where MALA is not: this applies, in particular, to one-dimensional super-Gaussian targets under the stated radial-growth assumptions.

We also derived high-dimensional optimal-scaling limits for MAIMLA on identically distributed and scaled independent product targets. The numerical investigations complement these results. In the Gaussian-mixture experiment, MAIMLA's sampling performance remains stable as the dimension increases, while MALA's effective sample size per iteration deteriorates. In the image-denoising experiment, both methods yield comparable reconstruction quality, while MAIMLA achieves higher effective sample sizes.

The Gaussian-mixture and image-denoising examples were chosen to illustrate MAIMLA's theoretical benefits over MALA, rather than to benchmark the cost of a general high-dimensional implicit solve. In the Gaussian-mixture experiment, the nonlinearity is confined to a single coordinate, reducing the implicit equation to a scalar solve; in the image-denoising experiment, the midpoint equation has a closed-form solution. For general targets, comparisons must account for the cost of solving the implicit step. In subsequent work, we will show that for specific targets arising in applications, such as change-of-measure problems from a Gaussian reference, MAIMLA can be more efficient than MALA even after this cost is taken into account.

\section*{Statement on AI use}
In this project AI was only used to aid in preparing the text of the manuscript. All AI assisted edits were reviewed by the authors who take full responsibility for the final content of the manuscript.

\appendix
\section{Midpoint Equation} \label{sec:midpoint}

For the proof of geometric ergodicity of IMLA and MAIMLA we used bounds on the midpoint $x_{\frac{1}{2}}=\frac{x+x^\star}{2}$, where $x^\star$ is the proposal of IMLA started from $x$. The midpoint can be seen as the solution of
\begin{equation}\label{eq:midpoint_equation}
    x_{\frac{1}{2}}+\frac{\delta}{2}\nabla U(x_{\frac{1}{2}}) = x+\sqrt{\frac{\delta}{2}}\xi.
\end{equation}

Note that if $\nabla U$ grows superlinearly then for large $x$ the behaviour of \eqref{eq:midpoint_equation} is dominated by
\begin{equation*}
    \frac{\delta}{2}\nabla U(x_{\frac{1}{2}}) \approx x+\sqrt{\frac{\delta}{2}}\xi.
\end{equation*}
Hence if $\nabla U$ grows like $\norm{x}^p$ then $\norm{x_{\frac{1}{2}}}$ is proportional to $\norm{x}^{\frac{1}{p}}$ for $x$ large. This is established in the following lemma.

\begin{lemma}\label{lem:midpoint_bounds}
    Suppose that $U\in \mathcal{C}^2$, let $x_{\frac{1}{2}}$ be the solution of \eqref{eq:midpoint_equation}, and set $z=x+\sqrt{\frac{\delta}{2}}\xi$.
    \begin{enumerate}
        \item \label{item:midpoint-lower-bound} Suppose that
        \begin{equation*}
            \norm{\nabla U(x)} \leq c_p\norm{x}^p+C_p.
        \end{equation*}
        Then
        \begin{equation}\label{eq:midpoint_lower_bound}
            \norm{x_{\frac{1}{2}}}
            \geq \min\left\{\frac{1}{2}\norm{z},
            \left(\left(\frac{1}{\delta c_p}\norm{z}-\frac{C_p}{c_p}\right)^+\right)^{\frac{1}{p}}\right\}.
        \end{equation}
        Here $\left(\cdot\right)^+ = \max\{\cdot,0\}$ denotes the positive part.
        \item \label{item:midpoint-upper-bound} Suppose that
        \begin{equation*}
            \norm{\nabla U(x)} \geq C_1\norm{x}^p+c_1.
        \end{equation*}
        Then
        \begin{equation*}
            \norm{x_{\frac{1}{2}}}
            \leq \max\left\{\left(\frac{4}{\delta C_1}\right)^{\frac{1}{p-1}},
            \left(\frac{4}{\delta C_1}(\norm{z}-\frac{\delta}{2}c_1)\right)^{\frac{1}{p}}\right\}.
        \end{equation*}
    \end{enumerate}
\end{lemma}

\begin{proof}
\textit{Proof of part~\ref{item:midpoint-lower-bound}.}
Using \eqref{eq:midpoint_equation}, we have
\begin{equation*}
    \norm{z}
    \leq \norm{x_{\frac{1}{2}}}
    + \frac{\delta}{2}\norm{\nabla U(x_{\frac{1}{2}})}.
\end{equation*}
If $\norm{x_{\frac{1}{2}}} \leq \frac{1}{2}\norm{z}$, then
\begin{equation*}
    \norm{z}
    \leq \frac{1}{2}\norm{z}
    + \frac{\delta}{2}\norm{\nabla U(x_{\frac{1}{2}})},
\end{equation*}
which rearranges to
\begin{equation*}
    \norm{z} \leq \delta \norm{\nabla U(x_{\frac{1}{2}})}.
\end{equation*}
Using the growth bound on $\nabla U$, we obtain
\begin{equation*}
    \norm{z}
    \leq \delta \norm{\nabla U(x_{\frac{1}{2}})}
    \leq \delta C_p + c_p\delta \norm{x_{\frac{1}{2}}}^p,
\end{equation*}
which rearranges to
\begin{equation*}
    \norm{x_{\frac{1}{2}}}
    \geq \left(\frac{1}{\delta c_p}\norm{z}-\frac{C_p}{c_p}\right)^{\frac{1}{p}}.
\end{equation*}
The right-hand side is positive only when $\norm{z}$ is sufficiently large; otherwise, we bound it by zero. Therefore, in general,
\begin{equation*}
    \norm{x_{\frac{1}{2}}}
    \geq \min\left\{\frac{1}{2}\norm{z},
    \left(\left(\frac{1}{\delta c_p}\norm{z}-\frac{C_p}{c_p}\right)^+\right)^{\frac{1}{p}}\right\}.
\end{equation*}

\textit{Proof of part~\ref{item:midpoint-upper-bound}.}
Using the bound and \eqref{eq:midpoint_equation}, we have
\begin{equation}\label{eq:upper_bound_p_power}
    \frac{\delta}{2} c_1+\frac{\delta}{2}C_1\norm{x_{\frac{1}{2}}}^p
    \leq \frac{\delta}{2}\norm{\nabla U(x_{\frac{1}{2}})}
    = \norm{z-x_{\frac{1}{2}}}
    \leq \norm{z} +\norm{x_{\frac{1}{2}}}.
\end{equation}

If $\norm{x_{\frac{1}{2}}} \leq \left(\frac{4}{\delta C_1}\right)^{\frac{1}{p-1}}$, then the required bound holds. We may therefore assume that
\begin{equation*}
    \norm{x_{\frac{1}{2}}} \geq \left(\frac{4}{\delta C_1}\right)^{\frac{1}{p-1}},
\end{equation*}
which implies
\begin{equation*}
    \frac{\delta C_1}{4}\norm{x_{\frac{1}{2}}}^p \geq \norm{x_{\frac{1}{2}}}.
\end{equation*}
Applying this inequality to \eqref{eq:upper_bound_p_power} gives
\begin{equation*}
    \frac{\delta}{2}c_1+\frac{\delta}{2}C_1\norm{x_{\frac{1}{2}}}^p
    \leq \norm{z} +\norm{x_{\frac{1}{2}}}
    \leq \norm{z} +\frac{\delta C_1}{4}\norm{x_{\frac{1}{2}}}^p.
\end{equation*}
After rearranging, we obtain
\begin{equation*}
    \frac{\delta}{4}C_1\norm{x_{\frac{1}{2}}}^p \leq \norm{z} -\frac{\delta}{2}c_1.
\end{equation*}
Hence,
\begin{equation*}
    \norm{x_{\frac{1}{2}}}
    \leq \left(\frac{4}{\delta C_1}(\norm{z}-\frac{\delta}{2}c_1)\right)^{\frac{1}{p}},
\end{equation*}
and the desired bound follows.
\end{proof}

\begin{lemma}\label{lem:radial_inner_bound}
    Suppose $U\in \mathcal{C}^2$ is of the form \eqref{eq:radial_condition}. Then
    \begin{align*}
        x^{\top}x_{\frac{1}{2}}
        &\geq \frac{\delta}{2}C(1+p)\norm{x_{\frac{1}{2}}}^{p+1}
        + \norm{x_{\frac{1}{2}}}^2
        - \delta K_p\norm{x_{\frac{1}{2}}}^p
        - K_p\left(1+\norm{\xi}^{\frac{p}{p-1}}\right),
    \end{align*}
    for a constant $K_p>0$ depending only on $p$, $C_0$, and $c_0$.
\end{lemma}

\begin{proof}
    To analyse this quantity, we multiply \eqref{eq:midpoint_equation} by $x_{\frac{1}{2}}^{\top}$ to get
    \begin{equation*}
        \norm{x_{\frac{1}{2}}}^2
        +\frac{\delta}{2}x_{\frac{1}{2}}^{\top}\nabla U(x_{\frac{1}{2}})
        = x_{\frac{1}{2}}^{\top}x
        +\sqrt{\frac{\delta}{2}}x_{\frac{1}{2}}^{\top}\xi.
    \end{equation*}
    By \eqref{eq:radial_condition}, we have
    \begin{equation*}
        \nabla U(x) = C(1+p)\norm{x}^{p-1}x + \nabla U_0(x),
    \end{equation*}
    which gives the following expression for $x^{\top}x_{\frac{1}{2}}$:
    \begin{equation*}
        x^{\top}x_{\frac{1}{2}}
        = \frac{\delta}{2}C(1+p)\norm{x_{\frac{1}{2}}}^{p+1}
        + \norm{x_{\frac{1}{2}}}^2
        + \frac{\delta}{2}x_{\frac{1}{2}}^{\top}\nabla U_0(x_{\frac{1}{2}})
        - \sqrt{\frac{\delta}{2}}x_{\frac{1}{2}}^{\top}\xi.
    \end{equation*}
    By Young's inequality,
    \begin{equation*}
        x^{\top}x_{\frac{1}{2}}
        \geq \frac{\delta}{2}C(1+p)\norm{x_{\frac{1}{2}}}^{p+1}
        + \norm{x_{\frac{1}{2}}}^2
        - \frac{\delta}{2p}\norm{x_{\frac{1}{2}}}^p
        - \frac{\delta(p-1)}{2p}\norm{\nabla U_0(x_{\frac{1}{2}})}^{\frac{p}{p-1}}
        - \frac{1}{p}\left(\frac{\delta}{2}\right)^{\frac{p}{2}}\norm{x_{\frac{1}{2}}}^p
        - \frac{p-1}{p}\norm{\xi}^{\frac{p}{p-1}}.
    \end{equation*}
    By using the assumption on $\nabla U_0$,
    \begin{align*}
        x^{\top}x_{\frac{1}{2}}
        &\geq \frac{\delta}{2}C(1+p)\norm{x_{\frac{1}{2}}}^{p+1}
        + \norm{x_{\frac{1}{2}}}^2
        - \frac{\delta}{2p}\norm{x_{\frac{1}{2}}}^p
        - \frac{\delta(p-1)}{2p}C_0^{\frac{p}{p-1}}\norm{x_{\frac{1}{2}}}^{p}
        - \frac{\delta(p-1)}{2p}c_0^{\frac{p}{p-1}}\\
        &- \frac{1}{p}\left(\frac{\delta}{2}\right)^{\frac{p}{2}}\norm{x_{\frac{1}{2}}}^p
        - \frac{p-1}{p}\norm{\xi}^{\frac{p}{p-1}}.
    \end{align*}
    Collecting terms, this gives
    \begin{align*}
        x^{\top}x_{\frac{1}{2}}
        &\geq \frac{\delta}{2}C(1+p)\norm{x_{\frac{1}{2}}}^{p+1}
        + \norm{x_{\frac{1}{2}}}^2
        - \delta K_p\norm{x_{\frac{1}{2}}}^p
        - K_p\left(1+\norm{\xi}^{\frac{p}{p-1}}\right),
    \end{align*}
    for a constant $K_p>0$ depending only on $p$, $C_0$, and $c_0$.
\end{proof}

\end{document}